# ASYMPTOTIC INFERENCE FOR SEMIPARAMETRIC ASSOCIATION MODELS

BY GERHARD OSIUS

*Universität Bremen*

Association models for a pair of random elements $X$ and $Y$ (e.g., vectors) are considered which specify the odds ratio function up to an unknown parameter $\boldsymbol{\theta}$. These models are shown to be semiparametric in the sense that they do not restrict the marginal distributions of $X$ and $Y$. Inference for the odds ratio parameter $\boldsymbol{\theta}$ may be obtained from sampling either $Y$ conditionally on $X$ or vice versa. Generalizing results from Prentice and Pyke, Weinberg and Wacholder and Scott and Wild, we show that asymptotic inference for $\boldsymbol{\theta}$ under sampling conditional on $Y$ is the same as if sampling had been conditional on $X$. Common regression models, for example, generalized linear models with canonical link or multivariate linear, respectively, logistic models, are association models where the regression parameter $\boldsymbol{\beta}$ is closely related to the odds ratio parameter $\boldsymbol{\theta}$. Hence inference for $\boldsymbol{\beta}$ may be drawn from samples conditional on $Y$ using an association model.

**1. Introduction and outline.** A common approach to describe the relationship between a random output variable $Y$ of interest (e.g., a health status) and a random input vector $X$ (e.g., consumption of tobacco, alcohol and other risk factors) is by means of a parametric regression model which specifies the conditional distribution of $Y$ given $X = x$ up to an unknown parameter vector. In the most simple case $Y$ is an indicator (e.g., for the presence of a disease) and the conditional distribution is binomial $B(1, p(x))$. The popular logistic regression model relates the logistic transform of $p(x)$ and a vector $\mathbf{z} = h(x) \in \mathbb{R}^S$ of covariates—obtained from $x$ by a suitable function $h$—through logit $p(x) = \gamma + \mathbf{z}^T \boldsymbol{\theta}$ with parameters $\gamma \in \mathbb{R}$ and $\boldsymbol{\theta} \in \mathbb{R}^S$. The appropriate sampling scheme for this model is to sample $Y$ conditionally on $X = x$ for specified values of $x$. In epidemiology this is called a *cohort study*, each of the $J$ cohorts being determined by its value $x$. In contrast, the so-called *case-control studies* are obtained by sampling $X$ conditional on $Y = 1$ (cases), respectively, $Y = 0$ (controls). An important result by Prentice and Pyke [12] briefly states that asymptotic inference for the parameter $\boldsymbol{\theta}$ (but *not* for $\gamma$) in a case-control study may be obtained as if the data came from a cohort study. Actually their work covers the multivariate logistic regression model (cf. Example 3) for a random variable $Y$ taking values in $\{0, 1, \ldots, K\}$ and was generalized by Scott and Wild [14] to multiplicative intercept models. Our aim is to extend these results to semiparametric *odds ratio models* (introduced in [9]) for random elements, including in particular random vectors $Y$ and $X$, each with continuous and/or discrete components. The odds ratio function $OR(x, y)$ for the joint density $p(x, y)$ of $X$ and $Y$ is defined as a cross-product ratio with respect to fixed reference values $x^{\circ}$ and $y^{\circ}$:

$$OR(x, y) = \frac{p(x, y) \cdot p(x^{\circ}, y^{\circ})}{p(x, y^{\circ}) \cdot p(x^{\circ}, y)}.$$

An equivalent description is given by the corresponding ratio for the *conditional* density $p(y \mid X = x)$ of $Y$ given $X$—or vice versa. Under mild assumptions the joint distribution of $(X, Y)$ is *uniquely* determined by the odds ratio function and the marginal distributions of $X$ and $Y$; compare [9] or [10]. And conversely, for any pair of marginal distributions for $X$ and $Y$ and an odds ratio function there exists a joint distribution having these properties. The odds ratio function thus captures the *complete* association structure of $X$ and $Y$ by ignoring the information contained in the marginal distributions. A parametric odds ratio model specifies only the odds ratio function up to an unknown parameter vector $\boldsymbol{\theta}$, that is,

$$\log OR(x, y) = \psi_{\boldsymbol{\theta}}(x, y).$$

This model is semiparametric in the sense that it does not restrict the marginal distributions of $X$ and $Y$, but only the association structure. An important class are *log-bilinear* association models where the log-odds ratio function is bilinear with respect to given transformations $\mathbf{z} = h_X(x)$ and $\mathbf{v} = h_Y(y)$, that is,

(1.1) $$\log OR(x, y) = \mathbf{z}^T \boldsymbol{\theta} \mathbf{v}.$$

In fact, some widely used regression models, for example, generalized linear models with canonical link function and multivariate linear, respectively, logistic regression models, have a log-bilinear association structure. The assumptions concerning the conditional distribution of $Y$ given $X$ in these regression models may be removed by passing to the corresponding log-bilinear odds ratio model. One advantage of odds ratio models over regression models is that inference about the odds ratio parameter $\boldsymbol{\theta}$ may be obtained from sampling $X$ conditionally on $Y$ or vice versa. To prove this, we first observe that maximum likelihood estimation is invariant under both conditional sampling schemes, that is, the estimate $\hat{\boldsymbol{\theta}}$ maximizing the conditional likelihood $L_{X|Y}$ for samples of $X$ given $Y$ also maximizes the corresponding conditional likelihood $L_{Y|X}$ for samples of $Y$ given $X$—and conversely. Generalizing the result in Prentice and Pike [12] and Scott and Wild [14], we show that the estimated asymptotic covariance matrix for $\hat{\boldsymbol{\theta}}$ is invariant under both conditional sampling schemes, too. Hence asymptotic inference concerning the odds ratio parameter $\boldsymbol{\theta}$ may be obtained from a sample drawn conditionally on $Y$ as if the sample had been drawn conditionally on $X$.

The paper is organized as follows. In Section 2 we establish that the joint distribution of $(X,Y)$ is uniquely determined by its odds ratio function and the marginal distributions (uniqueness theorem), and that each of these three components can vary independently of another (existence theorem). The latter result will be proved here under weaker assumptions than in [9] using a different approach. Association models are introduced in Section 3 and some widely used regression models are recognized having a log-bilinear association. Although log-bilinear association is a natural and common choice, we derive the main results for more general odds ratio models determined by

$$(1.2) \qquad \log OR(x,y) = G(\mathbf{z}, \mathbf{v}, \boldsymbol{\theta}),$$

where $G$ is a given (sufficiently smooth) function. Section 4 establishes that the maximum likelihood estimate $\hat{\boldsymbol{\theta}}$ is invariant under the usual sampling schemes: unconditional or conditional on $X$, respectively, $Y$. For log-bilinear association models the likelihood to maximize corresponds to a log-linear model for a suitable contingency table. Hence results on the existence and uniqueness as well as techniques to compute the estimate are already available.

Knowing that the estimate $\hat{\boldsymbol{\theta}}$ is invariant under conditional sampling given either $X$ or $Y$, we establish in several steps our main result, that its estimated asymptotic normal distribution is invariant, too. In Section 5 we consider sampling $X$ conditional on $Y$ but maximize the "reverse" conditional log-likelihood $\ell(\boldsymbol{\lambda})$—arising from conditioning $Y$ on $X$—with respect to $\boldsymbol{\lambda} = (\boldsymbol{\theta}, \boldsymbol{\gamma}^*)$, where $\boldsymbol{\gamma}^*$ is a nuisance parameter vector. For the information matrix $\mathbf{I}(\boldsymbol{\lambda}) = E(-D^2_{\boldsymbol{\lambda}\boldsymbol{\lambda}}\ell(\boldsymbol{\lambda}))$ we show that the submatrix $[\mathbf{I}^{-1}(\boldsymbol{\lambda})]_{\boldsymbol{\theta}\boldsymbol{\theta}}$ of $\mathbf{I}^{-1}(\boldsymbol{\lambda})$ corresponding to $\boldsymbol{\theta}$ is indeed the asymptotic covariance matrix of $\hat{\boldsymbol{\theta}}$. To establish the asymptotic normality of the estimate $\hat{\boldsymbol{\lambda}}$, we first prove its consistency in Section 6. Our asymptotic approach applies to a *fixed* set $\{y_0, \ldots, y_K\}$ of values for $Y$ to be conditioned upon and independent samples of size $n_k$ from each conditional distribution of $X$ given $Y = y_k$, such that $n = \sum_k n_k$ tends to infinity while the ratios $n_k/n$ remain *fixed*.

In Section 7 the asymptotic normality is derived more generally for *any* (weakly) consistent estimate $\hat{\boldsymbol{\lambda}}$ which solves the estimating equation at least approximately, that is, $D_{\boldsymbol{\lambda}}\ell(\hat{\boldsymbol{\lambda}}) = o_P(\sqrt{n})$. Using the observed information $\mathbf{J}(\hat{\boldsymbol{\lambda}}) = -D^2_{\boldsymbol{\lambda}\boldsymbol{\lambda}}\ell(\hat{\boldsymbol{\lambda}})$ as a consistent estimate of $\mathbf{I}(\boldsymbol{\lambda})$, we finally obtain the asymptotic normality of the odds ratio estimate

$$\hat{\boldsymbol{\theta}} \underset{\text{as.}}{\sim} N(\boldsymbol{\theta}, [\mathbf{J}^{-1}(\hat{\boldsymbol{\lambda}})]_{\boldsymbol{\theta}\boldsymbol{\theta}}).$$

The estimated asymptotic covariance matrix here is exactly the same as if sampling had been conditional on $X$ for the observed $x$-values.

We do not attempt to derive our results under the weakest possible assumptions but prefer a few easily interpretable conditions which will be verified for a log-bilinear association model under mild distributional assumptions. The approach adopted here is *symmetric* in $X$ and $Y$ so that interchanging $X$ with $Y$ in any argument entails its *dual*.

**2. The odds ratio function.** Consider arbitrary nonempty spaces $\Omega_X$, respectively $\Omega_Y$, with $\sigma$-algebras $\mathcal{B}_X$, respectively $\mathcal{B}_Y$, and denote the product $\sigma$-algebra on $\Omega = \Omega_X \times \Omega_Y$ by $\mathcal{B}$. Let $\mathcal{P}$ the space of all probability measures $P$ on $(\Omega, \mathcal{B})$ and denote the marginal distributions of $P$ on $\Omega_X$, respectively $\Omega_Y$, by $P^X$, respectively $P^Y$. The definition of an odds ratio function for $P$ requires a *positive* density with respect to a product measure and a natural choice is the product $P^{XY} = P^X \times P^Y$ of the marginals. This leads to the subspace of probability measures $P$ having a positive density with respect to $P^{XY}$, or equivalently, are dominated by and dominate $P^{XY}$:

$$\mathcal{P}_{\ll} = \left\{ P \in \mathcal{P} \,\Big|\, \frac{dP}{dP^{XY}} > 0 \right\} = \{ P \in \mathcal{P} \mid P \ll P^{XY} \ll P \}.$$

For any $P \in \mathcal{P}_{\ll}$ with density $p = dP/dP^{XY}$ its *odds ratio function* $OR_p$ with respect to fixed reference values $x^\circ \in \Omega_X$ and $y^\circ \in \Omega_Y$ is defined on $\Omega \times \Omega$ by

$$(2.1) \qquad OR_p(x,y) = \frac{p(x,y) \cdot p(x^\circ, y^\circ)}{p(x, y^\circ) \cdot p(x^\circ, y)}.$$

The choice of the dominating product measure $P^{XY}$ is not essential (cf. [9]): replacing $p$ by a *positive* density $p_\nu$ with respect to a product $\nu = \nu_X \times \nu_Y$ of $\sigma$-finite measures yields the same ratio (2.1). Since the density $p$ of $P$ is only unique up to almost sure equality, the same holds for the odds ratio function $OR_p$ of $P$, which nevertheless will also be denoted simply by $OR(P)$. The log-odds ratio function may be written in terms of the log-density

$$(2.2) \quad \log OR_p(x,y) = \log p(x,y) + \log p(x^\circ, y^\circ) - \log p(x, y^\circ) - \log p(x^\circ, y).$$

It is convenient to view any $P \in \mathcal{P}$ as a joint distribution of a pair $(X, Y)$ of random elements defined on some probability space with values in $\Omega$ and the odds ratio function of $(X, Y)$ is defined by $OR(X, Y) = OR(P)$.

To show that the odds ratio function completely characterizes the association between $X$ and $Y$, we have to restrict the joint distribution $P$ by requiring that its log-density $\log p$ is $P^{XY}$-integrable, or equivalently, that the *Kullback–Leibler information* [7]

$$I(P^{XY} \mid P) = \int \log\left(\frac{dP^{XY}}{dP}\right) dP^{XY}$$

is *finite*. Any $P$ in the subclass $\mathcal{P}_\int = \{P \in \mathcal{P}_\ll \mid I(P^{XY} \mid P) < \infty\}$ is uniquely determined by its marginal distributions and its odds ratio function.

THEOREM 1 (Uniqueness). *Any $P_1, P_2 \in \mathcal{P}_\int$ having the same marginals $P_1^X = P_2^X$, $P_1^Y = P_2^Y$ and the same odds ratio function $OR(P_1) = OR(P_2)$ agree: $P_1 = P_2$.*

For a proof one easily establishes $I(P_1 \mid P_2) = 0$ using (2.2); compare [10].

Next we want to "define" a distribution $P$ on $\Omega$ by specifying its marginal distributions and its (log) odds ratio function. For given distributions $\pi_X$ on $\Omega_X$ and $\pi_Y$ on $\Omega_Y$ and a measurable function $\psi$ on $\Omega$, we investigate under which conditions we can find a $P \in \mathcal{P}_\int$ with $P^X = \pi_X$, $P^Y = \pi_Y$ and $\log OR(P) = \psi$. First of all, $\psi$ has to satisfy the obvious constraints

CONDITION (OR1). $\psi(x, y^\circ) = 0$, $\psi(x^\circ, y) = 0$ for all $x, y$.

Furthermore from $P \in \mathcal{P}_\int$ and (2.2) we obtain two *necessary* integrability conditions:

CONDITION (E1). $\psi$ is $\pi_X \times \pi_Y$-integrable.

CONDITION (E2). There exists $\pi_X$-integrable $\beta: \Omega_X \longrightarrow \mathbb{R}$ $\pi_Y$-integrable $\gamma: \Omega_Y \longrightarrow \mathbb{R}$ functions such that $\exp(\psi - \beta - \gamma)$ is $\pi_X \times \pi_Y$-integrable.

These conditions are also sufficient for the existence of the wanted $P \in \mathcal{P}_\int$.

THEOREM 2 (Existence). *For distributions $\pi_X$ on $\Omega_X$ and $\pi_Y$ on $\Omega_Y$ and a measurable function $\psi$ on $\Omega \times \Omega$ the following statements are equivalent:*

(a) *There exists $P \in \mathcal{P}_\int$ with $P^X = \pi_X$, $P^Y = \pi_Y$ and $\log OR(P) = \psi$.*

(b) *There exists $P \in \mathcal{P}_f$ with $\log OR(P) = \psi$.*
(c) *$\psi$ satisfies Conditions (OR1), (E1) and (E2).*

The proof is given in Appendix A.1. A few remarks are in order.

1. Conditions (E1) and (E2) hold for *bounded* $\psi$, for example, for *continuous* $\psi$ and *compact* $\Omega$.
2. The integrability of $\exp(\psi - \beta - \gamma)$ in Condition (E2) holds if $\psi \leq \beta + \gamma$. And if even $|\psi| \leq \beta + \gamma$, then Condition (E1) follows, too.
3. For *finite* $\Omega_Y$ (*or* $\Omega_X$) Condition (E1) implies Condition (E2) for $\beta(x) = \sum_y |\psi(x,y)|$ and $\gamma = 0$.
4. Although $P$ is uniquely determined by Theorem 1, there is no *explicit* formula for $P$ available. In the proof $P$ is given by an $I$-projection, which can only be obtained as a limit in an iterative procedure. Only for *binary* $Y$ (and vector-valued $X$) the distribution $P$ is easily available; compare [1] or [9].
5. A stronger version of Condition (E2) requiring $\exp(\psi - \beta)$ and $\exp(\psi - \gamma)$ to be integrable was used in [9, 10] to obtain $P$ as a limit of an iterative proportional fitting procedure.
6. For *finite* spaces $\Omega_X$ and $\Omega_Y$ this result has long been known; compare [11], Section 3.4.

**3. Association models.** An association model for the joint distribution $P$ of $(X,Y)$ only restricts the odds ratio function of $P$ and leaves the marginal distributions of $X$ and $Y$ arbitrary. To formulate such a model we assume that $P$ has a *positive* density with respect to a fixed product measure $\nu = \nu_X \times \nu_Y$ of $\sigma$-finite measures $\nu_X$, respectively $\nu_Y$, on $\Omega_X$, respectively $\Omega_Y$. Hence $P$ is restricted to the class $\mathcal{P}^{XY} = \{P \in \mathcal{P} \mid P \ll \nu \ll P\} \subset \mathcal{P}_{\ll}$, which also restricts the marginal distribution $P^X$ of $X$ to

$$\mathcal{P}^X = \{\pi_X \text{ probability measure on } \Omega_X \mid \pi_X \ll \nu_X \ll \pi_X\},$$

and the marginal $P^Y$ to the corresponding $\mathcal{P}^Y$. From now on all densities on $\Omega$, respectively $\Omega_X, \Omega_Y$ are taken with respect to the dominating measure $\nu$, respectively $\nu_X, \nu_Y$.

We consider parametric association models indexed by a parameter vector $\boldsymbol{\theta} \in \mathbb{R}^S$. For any $\boldsymbol{\theta}$ let $\psi_{\boldsymbol{\theta}}$ be a measurable function on $\Omega$ satisfying Condition (OR1). The parametric odds ratio model restricts the log-odds ratio function of $P$ to $\log OR(P) = \psi_{\boldsymbol{\theta}}$ for some $\boldsymbol{\theta}$. To guarantee for any $\boldsymbol{\theta}$ and any marginals $\pi_X, \pi_Y$ the existence of a joint distribution $P$ with $\psi_{\boldsymbol{\theta}} = \log OR(P)$ and these marginals, we assume the following bounding condition:

CONDITION (OR2). *There exist nonnegative measurable functions $\tilde{\psi}_X$ on $\Omega_X$ and $\tilde{\psi}_Y$ on $\Omega_Y$ with $|\psi_{\boldsymbol{\theta}}(x,y)| \leq [\tilde{\psi}_X(x) + \tilde{\psi}_Y(y)] \cdot \|\boldsymbol{\theta}\|$ for all $\boldsymbol{\theta}, x, y$.*

Furthermore we restrict $\pi_X$ to the class $\mathcal{P}_f^X = \{\pi_X \in \mathcal{P}^X \mid \tilde{\psi}_X \text{ is } \pi_X\text{-integrable}\}$ and $\pi_Y$ to the corresponding class $\mathcal{P}_f^Y$. Condition (c) in Theorem 2 holds for any $\pi_X \in \mathcal{P}_f^X$, $\pi_Y \in \mathcal{P}_f^Y$ and $\boldsymbol{\theta}$, and hence there exists a unique $P \in \mathcal{P}_f$ with $P^X = \pi_X, P^Y = \pi_Y$ and $\log OR(P) = \psi_{\boldsymbol{\theta}}$. Thus a *parametric association model* (PAM) for distributions $P$ in $\mathcal{P}_f^{XY} = \mathcal{P}^{XY} \cap \mathcal{P}_f$ is specified by the requirements

$$(3.1) \qquad \log OR(P) \in \{\psi_{\boldsymbol{\theta}} \mid \boldsymbol{\theta} \in \mathbb{R}^S\}, \qquad P^X \in \mathcal{P}_f^X, \ P^Y \in \mathcal{P}_f^Y.$$

This is a semiparametric model for the joint distribution $P$ since the marginals are only slightly restricted by integrability conditions. By (2.2) a density $p(x,y)$ of $P \in \mathcal{P}_f^{XY}$ satisfying (3.1) can be parametrized as

$$(3.2) \qquad \log p(x,y) = \alpha + \beta(x) + \gamma(y) + \psi_{\boldsymbol{\theta}}(x,y)$$

with $\alpha \in \mathbb{R}$ and integrable functions $\beta$ and $\gamma$. Identifiability may be achieved through the constraints $\beta(x^\circ) = 0$ and $\gamma(y^\circ) = 0$, which will be assumed here. The integration constant $\alpha$ is determined by

$$\alpha = -\log \int \exp(\beta + \gamma + \psi_{\boldsymbol{\theta}}) \, d\nu$$

and marginal density $p^X(x)$ of $P^X$ is given by

$$\log p^X(x) = \alpha + \beta(x) + \delta(x), \qquad \delta(x) = \log \left[ \int \exp(\gamma(y) + \psi_{\boldsymbol{\theta}}(x,y)) \, d\nu_Y(y) \right].$$

The conditional distribution of $Y$ given $X = x$ belongs to $\mathcal{P}^Y$ and the conditional density $p(y \mid X = x)$ satisfies

$$(3.3) \qquad \log p(y \mid X = x) = \gamma(y) + \psi_{\boldsymbol{\theta}}(x,y) - \delta(x).$$

The integration constant $\delta(x)$ can be removed by passing to the density ratio

$$(3.4) \qquad \log \frac{p(y \mid X = x)}{p(y^\circ \mid X = x)} = \gamma(y) + \psi_{\boldsymbol{\theta}}(x,y).$$

Equation (3.4) may be viewed as a "regression model." Conversely, suppose a model for $P$ is specified by (3.4) with an *arbitrary* integrable function $\gamma$ and the parametric family $\psi_{\boldsymbol{\theta}}$. Then $\log OR(P) = \psi_{\boldsymbol{\theta}}$ and hence the model (3.4) is semiparametric in the sense that it does not restrict the marginal distributions $P^X$ and $P^Y$—provided they belong to the class $\mathcal{P}_f^X$, respectively $\mathcal{P}_f^Y$. In the latter case the regression model (3.4) is in fact *equivalent* to the association model (3.1). Note that for finite $\Omega_Y$ and counting measure $\nu_Y$ the integrability condition imposed by $P^Y \in \mathcal{P}_f^Y$ always holds.

An important class of parametric association models are *log-bilinear association (LBA) models* with respect to measurable maps $h_X : \Omega_X \longrightarrow \mathbb{R}^{K_X}$ and $h_Y : \Omega_Y \longrightarrow \mathbb{R}^{K_Y}$, which will always be chosen here such that $h_X(x^\circ) = \mathbf{0}$ and $h_Y(y^\circ) = \mathbf{0}$. The parameter $\boldsymbol{\theta}$ is a $K_X \times K_Y$-matrix and the log-odds ratio function is bilinear in the transformed variables $h_X(x)$ and $h_Y(y)$

$$(3.5) \qquad \psi_{\boldsymbol{\theta}}(x,y) = h_X(x)^T \boldsymbol{\theta} h_Y(y) \qquad \text{for all } x, y.$$

Since $|h_X(x)^T \boldsymbol{\theta} h_Y(y)| \leq \|h_X(x)\| \cdot \|h_Y(y)\| \cdot \|\boldsymbol{\theta}\|$, Condition (OR2) holds for $\tilde{\psi}_X(x) = \|h_X(x)\|^2$ and $\tilde{\psi}_Y(y) = \|h_Y(y)\|^2$. And the integrability condition in $\mathcal{P}_f^X$ and $\mathcal{P}_f^Y$ states that the second moments $E(\|h_X(X)\|^2)$ and $E(\|h_Y(Y)\|^2)$ are finite. Any submodel of (3.5) specified by a linear restriction of the form $\boldsymbol{\theta} = \mathbf{A}^T \boldsymbol{\theta}^* \mathbf{B}$ with given matrices $\mathbf{A}$, $\mathbf{B}$ and parameter matrix $\boldsymbol{\theta}^*$ yields a log-bilinear association too, with respect to $h_X^* = \mathbf{A} h_X$, $h_Y^* = \mathbf{B} h_Y$.

Association models have been introduced long ago in the context of contingency tables, that is, when both $X$ and $Y$ have a *finite* range; see [4] for a review. The *"RC association models"* and *"RC correlation models"* in [4] are both association models in our sense, the former (but not the latter) being log-bilinear. Extensions of these models to multivariate contingency tables studied in Gilula and Haberman [3] also satisfy (3.1). Goodman [4] has generalized the bivariate normal distribution to a bivariate log-bilinear model in our sense, but did not establish its semiparametric nature. Returning to our primary focus, namely *general random vectors* $X$ and $Y$, the following examples reveal that the association structure of some widely used regression models is in fact *log-bilinear*.

EXAMPLE 1 (Generalized linear models). Let $Y$ be a univariate random variable, $X$ an $R$-dimensional random vector and suppose that the conditional density of $Y$ given $X = x$ belongs to the exponential family $p(y \mid X = x) = \exp\{a(\phi)^{-1}[y \cdot \tau(x) - b(\tau(x))] + c(y, \phi)\}$ with suitable functions $a, b, c, \tau$ and a (dispersion) parameter $\phi$; compare [8]. Then the log-odds ratio function has the form $\psi(x, y) = a(\phi)^{-1} \cdot [\tau(x) - \tau(x^\circ)] \cdot [y - y^\circ]$ and $\tau(x)$ is a strictly monotone function of the conditional expectation $\mu(x) = E(Y \mid X = x)$, namely $\tau(x) = \lambda(\mu(x))$, where $\lambda^{-1} = b'$. A *generalized linear model* specifies the conditional expectation via a link function $g$:

$$(3.6) \qquad g(\mu(x)) = \alpha + \mathbf{z}^T \boldsymbol{\beta},$$

where $\mathbf{z} = h_X(x) \in \mathbb{R}^S$ is a known vector of *formal* covariates (obtained from $x$ by a given function $h_X$) and $\alpha \in \mathbb{R}$, $\boldsymbol{\beta} \in \mathbb{R}^S$ are unknown parameters. For $G = \lambda \circ g^{-1}$ and $h_X(x^\circ) = \mathbf{0}$ the log-odds ratio function is $\psi(x, y) = a(\phi)^{-1} \cdot [G(\alpha + \mathbf{z}^T \boldsymbol{\beta}) - G(\alpha)] \cdot [y - y^\circ]$. If the *canonical link* $g = \lambda^{-1}$ is chosen, then

$$(3.7) \qquad \psi(x, y) = \mathbf{z}^T \boldsymbol{\theta} [y - y^\circ]$$

is of the form (3.5) with $h_Y(y) = y - y^\circ$ and parameter $\boldsymbol{\theta} = a(\phi)^{-1}\boldsymbol{\beta}$. Note that the intercept $\alpha$ is no longer present in (3.7). Taking the log-bilinear association model (3.7) instead of (3.6) weakens the distributional assumption while still including the regression parameter $\boldsymbol{\beta}$ up to a positive constant $a(\phi)^{-1}$. In particular a linear hypothesis $\mathbf{C}\boldsymbol{\beta} = \mathbf{0}$ with a given matrix $\mathbf{C}$ is equivalent to $\mathbf{C}\boldsymbol{\theta} = \mathbf{0}$, and for a vector $\mathbf{c}$ a one-sided hypothesis $\mathbf{c}^T\boldsymbol{\beta} > 0$ is equivalent to $\mathbf{c}^T\boldsymbol{\theta} > 0$. Generalized linear models with canonical link are often used. First of all, normal conditional distributions $N(\mu(x), \sigma^2)$ of $Y$ yield the classical linear model with $a(\phi) = \sigma^2$. Second, binomial conditional distributions $B(\mu(x), 1)$ lead to logistic regression models. And finally, for Poisson conditional distributions $Pois(\mu(x))$ log-linear models are obtained. Note that for the latter two models we have $a(\phi) = 1$ and hence $\boldsymbol{\theta} = \boldsymbol{\beta}$.

The above semiparametric nature of the logistic regression model has been noticed before; compare Breslow, Robins and Wellner [1], who established its semiparametric efficiency under case-control sampling. However, the logistic regression model is the only one among generalized linear models for binary $Y$ which is equivalent to an association model (3.1); compare [9] or Example 2 below. And the resulting relation between the two conditional densities (given $X$, resp., $Y$) has been noticed before by Kagan [6].

EXAMPLE 2 (Multivariate linear logistic regression). Extending univariate logistic regression to the multivariate case, suppose $Y$ (e.g., a disease status) takes values in $\Omega_Y = \{0, 1, \ldots, K\}$, $K \geq 1$, and $X$ is an $R$-dimensional vector of *observed* covariates. Then $\mathcal{L}(Y \mid X = x)$ is a multinomial distribution $M_{K+1}(1, \pi(x))$ with $K + 1$ classes and probabilities $\pi_k(x) = P(Y = k \mid X = x) > 0$. Using the multivariate logistic transformation logit $\pi_k(x) = \log(\pi_k(x)/\pi_0(x))$ of $\pi(x)$, the *linear logistic regression model* is given by

$$(3.8) \qquad \text{logit } \pi_k(x) = \gamma_k + \mathbf{z}^T \boldsymbol{\theta}_k, \qquad k = 1, \ldots, K,$$

where $\mathbf{z} = h_X(x) \in \mathbb{R}^S$ is as above a vector of *formal* covariates and $\gamma_k \in \mathbb{R}$, $\boldsymbol{\theta}_k \in \mathbb{R}^S$ are unknown parameters. Choosing $y^\circ = 0$, the log-odds ratio function is

$$(3.9) \qquad \psi(x, y) = h_X(x)^T \boldsymbol{\theta}_k = h_X(x)^T \boldsymbol{\theta} h_Y(k),$$

where $\boldsymbol{\theta} = (\boldsymbol{\theta}_1, \ldots, \boldsymbol{\theta}_K)$ is an $S \times K$ parameter matrix, and the function $h_Y : \Omega_Y \longrightarrow \mathbb{R}^K$ maps $k > 0$ to the $k$th unit vector $\mathbf{e}_k$ and $h_Y(0) = \mathbf{0}$. Hence the linear logistic regression model is equivalent to the log-bilinear association model (3.9)—provided $E(\|h_X(X)\|^2)$ is finite. As mentioned above, this also holds for submodels given by linear constraints, for example, $\boldsymbol{\theta}_k = \boldsymbol{\theta}^*$ for all $k > 0$. Although the model (3.8) has been known for a long time, its *semiparametric* character (based on Theorem 2) does not seem to have been established before for $K > 2$.

Replacing $\mathbf{z}^T \boldsymbol{\theta}_k$ by an arbitrary function $g(\mathbf{z}, \boldsymbol{\theta}_k)$ leads to a *general logistic regression* model

$$\text{logit } \pi_k(x) = \gamma_k + g(\mathbf{z}, \boldsymbol{\theta}_k), \qquad k = 1, \ldots, K,$$

which is *equivalent* to the log-odds ratio model

$$\psi(x, y) = g(h_X(x), \boldsymbol{\theta}_k) = g(h_X(x), \boldsymbol{\theta} h_Y(k)).$$

EXAMPLE 3 (Multivariate linear regression). Let $Y$ and $X$ be random vectors taking values in $\mathbb{R}^K$, respectively $\mathbb{R}^R$, and suppose that the conditional distribution of $Y$ given $X$ is multivariate normal,

(3.10) $$\mathcal{L}(Y \mid X = x) = N_K(\mu_Y(x), \boldsymbol{\Sigma}),$$

such that the conditional covariance matrix $\boldsymbol{\Sigma}$ is nonsingular and does not depend on $x$. From the conditional log-density

$$\log p(y \mid X = x) = -\tfrac{1}{2}[\log[(2\pi)^K \det(\boldsymbol{\Sigma})] + [y - \mu_Y(x)]^T \boldsymbol{\Sigma}^{-1}[y - \mu_Y(x)]]$$

the log-odds ratio function with respect to $y^\circ = \mathbf{0}$ is $\psi(x, y) = [\mu_Y(x) - \mu_Y(x^\circ)]^T \boldsymbol{\Sigma}^{-1} y$. The *multivariate linear regression model*

(3.11) $$\mu_Y(x) = \boldsymbol{\alpha} + \boldsymbol{\beta}^T \mathbf{z}$$

with covariates $\mathbf{z} = h_X(x) \in \mathbb{R}^S$ and $S \times K$ parameter matrix $\boldsymbol{\beta}$ has a log-bilinear association

(3.12) $$\psi(x, y) = h_X(x)^T \boldsymbol{\theta} y$$

with parameter matrix $\boldsymbol{\theta} = \boldsymbol{\beta} \boldsymbol{\Sigma}^{-1}$—assuming $h_X(x^\circ) = \mathbf{0}$. Note that the regression parameter $\boldsymbol{\beta}$ may only be recovered from $\boldsymbol{\theta}$ if the covariance matrix $\boldsymbol{\Sigma}$ is known. However, any linear hypothesis $\mathbf{C}\boldsymbol{\beta} = \mathbf{0}$ is equivalent to the corresponding hypothesis $\mathbf{C}\boldsymbol{\theta} = \mathbf{0}$, and the latter may be tested using the semiparametric association model (3.12) instead of the regression model (3.11) with the distributional assumption (3.10). If instead of (3.10) we allow the conditional covariance matrix to depend on $x$, that is, $\mathcal{L}(Y \mid X = x) = N_K(\mu_Y(x), \boldsymbol{\Sigma}(x))$, then (3.11) leads to $\psi(x, y) = h_X(x)^T \boldsymbol{\beta} \boldsymbol{\Sigma}^{-1}(x) y$, which is not bilinear.

The above examples reveal that important regression models may be generalized to log-bilinear association models by ignoring the distributional assumption for the conditional distribution. Although log-bilinear association is a natural candidate, we also consider the more general *association model*

(3.13) $$\psi_{\boldsymbol{\theta}}(x, y) = G(h_X(x), h_Y(y), \boldsymbol{\theta}) \qquad \text{for all } x, y,$$

given by a fixed function $G$ with $G(\mathbf{0}, -, -) = G(-, \mathbf{0}, -) = 0$. We assume throughout that the function $G$ satisfies the following regularity condition (although some results also hold under weaker assumptions):

CONDITION (R1). $G(\mathbf{z}, \mathbf{v}, \boldsymbol{\theta})$ is thrice continuously differentiable with respect to $\boldsymbol{\theta}$ for all $\mathbf{z} \in h_X[\Omega_X]$, $\mathbf{v} \in h_Y[\Omega_Y]$ and the derivatives are continuous in $\mathbf{z}$ and $\mathbf{v}$.

Further properties of the functions $h_X, h_Y$ and $G$ will be assumed later in Conditions (R2″) and (MC).

**4. Estimation.** For a given data set $(x_i, y_i)$ with $i = 1, \ldots, n$ we want to estimate the association parameter $\boldsymbol{\theta}$ of the model (3.13) under unconditional sampling from the joint distribution of $(X, Y)$ and conditional sampling of $Y$ given $X$ or vice versa. Not surprisingly the maximum likelihood estimate $\hat{\boldsymbol{\theta}}$ under *any* of these three sampling schemes may be obtained as a solution of the *same* estimating equation.

4.1. *Unconditional sampling.* For unconditional sampling the data set $(x_i, y_i)$ is an independent sample from the joint distribution of $(X, Y)$. Suppose there are $J + 1 > 1$ different $x$-values and $K + 1 > 1$ different $y$-values observed and denote the corresponding subsets of $\Omega_X$ and $\Omega_Y$ by $\Omega_X^* = \{x_{(0)}, \ldots, x_{(J)}\}$ and $\Omega_Y^* = \{y_{(0)}, \ldots, y_{(K)}\}$. If $r_{jk}$ is the observed frequency of $(x_{(j)}, y_{(k)})$, then the likelihood is

$$L_{XY} = \prod_{j=0}^{J} \prod_{k=0}^{K} p(x_{(j)}, y_{(k)})^{r_{jk}} = L_{X|Y} \cdot L_Y$$

with a conditional and a marginal likelihood

$$(4.1) \quad L_{X|Y} = \prod_{k=0}^{K} \prod_{j=0}^{J} p(x_{(j)} \mid Y = y_{(k)})^{r_{jk}}, \qquad L_Y = \prod_{k=0}^{K} p^Y(y_{(k)})^{r_{+k}}$$

(the subscript "+" indicates summation over the replaced index). The model does not restrict the marginal distribution of $Y$ and hence the empirical density with respect to counting measure $\nu_Y^*$ on $\Omega_Y^*$,

$$(4.2) \qquad \hat{p}^Y(y_{(k)}) = \frac{1}{n} r_{+k} \qquad \text{for } k = 0, \ldots, K$$

is the usual nonparametric estimate. If we restrict the distribution $P^Y$ to the class $\mathcal{P}_Y^*$ of all distributions with finite support $\Omega_Y^*$, then $L_Y$ is a multinomial likelihood which attains its maximum for (4.2). Hence, for estimation purposes we may restrict the marginal $P^Y$ to $\mathcal{P}_Y^*$ and maximization of $L_{XY}$ is equivalent to separate maximization of $L_{X|Y}$ and $L_Y$, because the latter two have no common parameters.

Interchanging $X$ and $Y$, we split the likelihood as $L_{XY} = L_{Y|X} \cdot L_X$ and by the above argument we may *additionally* restrict $P^X$ to the class $\mathcal{P}_X^*$ of all distributions with finite support $\Omega_X^*$. Under these restrictions for *both* $P^X$ and $P^Y$ the likelihood $L_{XY}$ is a multinomial likelihood for the observed $(J+1) \times (K+1)$-contingency table $(r_{jk})$. Hence, estimation of $\boldsymbol{\theta}$ is reduced to a multinomial model whose probabilities $p_{jk} = p(x_{(j)}, y_{(k)})$ satisfy the log-odds ratio model

$$\log(p_{jk}p_{00}/p_{j0}p_{0k}) = \psi_{\boldsymbol{\theta}}(x_{(j)}, y_{(k)}) =: \psi_{jk}(\boldsymbol{\theta}) \qquad \text{for all } j \text{ and } k$$

with respect to the reference values $x^\circ = x_{(0)}$ and $y^\circ = y_{(0)}$. The parametrization (3.2) now involves only a *finite* number of parameters

$$(4.3) \quad \log p_{jk} = \beta_j + \gamma_k + \psi_{jk}(\boldsymbol{\theta}) - \log\left(\sum_j \sum_k \exp[\beta_j + \gamma_k + \psi_{jk}(\boldsymbol{\theta})]\right),$$

namely $\beta_j = \beta(x_{(j)})$, $\gamma_k = \gamma(y_{(k)})$ and $\boldsymbol{\theta}$ with $\beta_0 = \gamma_0 = 0$. Instead of maximizing $L_{XY}$, it is typically preferable to maximize either $L_{Y|X}$ or $L_{X|Y}$ using the parametrization of the conditional probabilities $p_{k|j} = p_{jk}/p_{j+}$ or $p_{j|k} = p_{jk}/p_{+k}$ given by (3.3) and its dual

$$\log p_{k|j} = \gamma_k + \psi_{jk}(\boldsymbol{\theta}) - \delta_j, \qquad \log p_{j|k} = \beta_j + \psi_{jk}(\boldsymbol{\theta}) - \varepsilon_k,$$

where the parameters $\delta_j$, respectively, $\varepsilon_k$ are determined by the remaining ones.

4.2. *Conditional sampling.* When sampling is conditional on values for $Y$ taken from $\Omega_Y^* = \{y_{(0)}, \ldots, y_{(K)}\}$, say, then the data set $(x_i, y_i)$ with $i = 1, \ldots, n$ is partitioned into $K+1$ independent subsamples given by the values of $y_i$, such that each subsample $(x_i)$ with $y_i = y_{(k)}$ is an independent sample from the conditional distribution $\mathcal{L}(X \mid Y = y_{(k)})$. Instead of maximizing the appropriate likelihood $L_{X|Y}$ we can equivalently maximize the unconditional likelihood $L_{XY}$ or even the "reverse" conditional likelihood $L_{Y|X}$. The latter is preferable from a computational point of view, when the nuisance parameters $\gamma_k$ are less than those of $L_{X|Y}$, that is, for $K < L$. A dual argument applies if sampling is conditional on values for $X$ taken from $\Omega_X^* = \{x_{(0)}, \ldots, x_{(J)}\}$.

4.3. *Log-bilinear association.* In the log-bilinear association model (3.5), the odds ratios may be written as $\psi_{jk}(\boldsymbol{\theta}) = \mathbf{z}_j^T \boldsymbol{\theta} \mathbf{v}_k$ with $\mathbf{z}_j = h_X(x_{(j)})$ and $\mathbf{v}_k = h_Y(y_{(k)})$, or in matrix notation

$$\boldsymbol{\psi}(\boldsymbol{\theta}) = \mathbf{Z}\boldsymbol{\theta}\mathbf{V}^T \in \mathbb{R}^{J \times K}, \qquad \mathbf{Z} = (\mathbf{z}_{jl}) \in \mathbb{R}^{J \times K_X}, \qquad \mathbf{V} = (v_{kl}) \in \mathbb{R}^{K \times K_Y}.$$

Then (4.3) reduces to a log-linear model for the probabilities $p_{jk}$,

$$(4.4) \qquad \log p_{jk} = \alpha + \beta_j + \gamma_k + \mathbf{z}_j^T \boldsymbol{\theta} \mathbf{v}_k$$

induced by the covariates $\mathbf{z}_j, \mathbf{v}_k$ and results by Haberman [5] on the existence and uniqueness of maximum likelihood estimates in log-linear models apply. In particular the estimate $\hat{\mathbf{p}} = (\hat{p}_{jk})$ is unique (if it exists) and hence the estimate $\hat{\boldsymbol{\theta}}$ is unique too, provided the parameter $\boldsymbol{\theta}$ is identifiable.

For sampling conditional on $Y$, the values $y_{(k)}$ should be chosen such that the rank condition holds:

CONDITION (RK). The $K_Y \times K$-matrix $\mathbf{V}^T = (\mathbf{v}_1, \ldots, \mathbf{v}_K)$ has rank $K_Y$.

This condition will be assumed whenever the log-bilinear association model is used. Then a convenient reparametrization is available:

$$(4.5) \qquad \psi_{jk}(\boldsymbol{\theta}) = \mathbf{z}_j^T \tilde{\boldsymbol{\theta}}_k, \qquad \tilde{\boldsymbol{\theta}}_k = \boldsymbol{\theta} \mathbf{v}_k \in \mathbb{R}^{K_X}$$

with a $K_X \times K$ parameter-matrix $\tilde{\boldsymbol{\theta}} = (\tilde{\boldsymbol{\theta}}_1, \ldots, \tilde{\boldsymbol{\theta}}_K) = \boldsymbol{\theta} \mathbf{V}^T$. The observed matrix $\mathbf{Z}$ of covariates will typically have rank $K_X$ and then $\boldsymbol{\theta}$, respectively, $\tilde{\boldsymbol{\theta}}$ is uniquely determined by $\boldsymbol{\psi}(\boldsymbol{\theta}) = \mathbf{Z}\boldsymbol{\theta}\mathbf{V}^T = \mathbf{Z}\tilde{\boldsymbol{\theta}}$ and hence identifiable. In general identifiability of $\boldsymbol{\theta}$ is guaranteed by Condition (C3) in Section 6.

**5. Conditional likelihood.** Although the maximum likelihood estimate $\hat{\boldsymbol{\theta}}$ of the association parameter $\boldsymbol{\theta}$ may be obtained by maximizing either of the two conditional likelihoods, the stochastic properties of the latter depend on the sampling scheme. Let us now consider sampling conditional on $Y$—which can be preferable from a practical point of view (even for regression models)—and derive properties of the "reverse" likelihood $L_{Y|X}$. The advantage of $L_{Y|X}$ over the appropriate likelihood $L_{X|Y}$ is that it usually has fewer nuisance parameters since $K$ is fixed by the sampling design whereas $J$ will typically increase with the number of observations—unless $\Omega_X$ is finite. An important example for *finite* $\Omega_Y$ are case-control studies (called *choice-based samples* in econometrics) for which asymptotic inference on $\boldsymbol{\theta}$ in the (general) logistic regression model may be obtained as if sampling had been conditional on $X$; compare [12] and [14]. We want to extend these results to *arbitrary* $Y$ (e.g., *vectors* with *continuous* and/or *discrete* components) and *association models*.

Instead of a data set $(x_i, y_i)$ we now consider the underlying random elements. It is convenient to represent the sample as a compound vector of random elements $\mathbf{X} = (X_{ki})$ indexed by $k = 0, \ldots, K$ and $i = 1, \ldots, n_k$. Omitting now the parentheses in $y_{(k)}$ and $x_{(j)}$, each $X_{ki}$ is distributed as $X_k \sim \mathcal{L}(X \mid Y = y_k)$. As above $r_{jk}$ denotes the frequency of $(x_j, y_k)$ in the sample $(x_{ki}, y_k)$ and the empirical distribution on $\Omega_Y^* = \{y_0, \ldots, y_K\}$ is given by the proportions $\bar{r}_k = n_k/n$, where $n = n_+$ is the total sample size. Replacing in $P$ the marginal distribution of $Y$ by the empirical distribution (4.2) yields a joint distribution $P^*$ on $\Omega_X \times \Omega_Y^*$ given by the density $p^*$ with respect to $\nu_X \times \nu_Y^*$:

$$p^*(x, y_k) = \bar{r}_k \cdot p(x \mid Y = y_k) \qquad \text{for all } x, k.$$

The marginal density of $Y$ under $P^*$ is $p^{*Y}(y_k) = \bar{r}_k$ and the marginal, respectively, conditional density for $X$ is

$$p^{*X}(x) = \sum_{k=0}^{K} \bar{r}_k \cdot p(x \mid Y = y_k), \qquad \text{respectively,}$$

(5.1)

$$p_k^*(x) := p^*(y_k \mid X = x) = \frac{\bar{r}_k \cdot p(x \mid Y = y_k)}{p^{*X}(x)}.$$

Equation (3.3) yields the parametrization $\log p_k^*(x) = \gamma_k^* + \psi_\theta(x, y_k) - \delta^*(x)$ with nuisance parameters $\gamma_k^* = \gamma^*(y_k)$ and $\delta^*(x) = \log[\sum_l \exp(\gamma_l^* + \psi_\theta(x, y_l))]$, hence

(5.2) $$p_k^*(x) = \frac{\exp[\gamma_k^* + \psi_\theta(x, y_k)]}{[\sum_l \exp[\gamma_l^* + \psi_\theta(x, y_l)]]}.$$

Choosing the reference value $y^\circ = y_0$ we have $\gamma_0^* = 0$, and the nuisance parameter is $\boldsymbol{\gamma}^* = (\gamma_1^*, \ldots, \gamma_K^*) \in \mathbb{R}^K$. Finally, the logarithm of the conditional likelihood $L_{Y|X}$ may be written in terms of the compound parameter vector $\boldsymbol{\lambda} := (\boldsymbol{\theta}, \boldsymbol{\gamma}^*) \in \mathbb{R}^{S+K}$:

$$\ell(\boldsymbol{\lambda}) := \log L_{Y|X} = \sum_{k=0}^{K} \sum_{i=1}^{n_k} \log p_k^*(X_{ki}) \qquad \text{with}$$

(5.3)

$$\log p_k^*(X_{ki}) = \gamma_k^* + \psi_\theta(X_{ki}, y_k) - \log\left[\sum_{l=0}^{K} \exp(\gamma_l^* + \psi_\theta(X_{ki}, y_l))\right].$$

Notice that $\ell(\boldsymbol{\lambda})$ is the log-likelihood of the multivariate logistic regression model

(5.4) $$\operatorname{logit} p_k^*(x) = \gamma_k^* + \psi_\theta(x, y_k), \qquad k = 1, \ldots, K,$$

which is *nonlinear* in general. The estimate $\hat{\boldsymbol{\lambda}}$ maximizing $\ell(\boldsymbol{\lambda})$ satisfies

(5.5) $$D_{\boldsymbol{\lambda}} \ell(\boldsymbol{\lambda}) = \sum_{k=0}^{K} \sum_{i=1}^{n_k} D_{\boldsymbol{\lambda}} \log p_k^*(X_{ki}) = \mathbf{0},$$

where $D_{\boldsymbol{\lambda}}$ denotes the differential operator with respect to $\boldsymbol{\lambda}$. The basic stochastic properties of the solution of the estimating equation (5.5) depend on the moments of the estimating function $D_{\boldsymbol{\lambda}} \ell(\boldsymbol{\lambda})$ and its derivative. The first important property (proved in Appendix A.2) is that its expectation is zero—which is not obvious since $\ell(\boldsymbol{\lambda})$ is *not* the log-likelihood for the underlying sampling:

$$(5.6) \qquad E[D_{\boldsymbol{\lambda}}\ell(\boldsymbol{\lambda})] = \sum_k n_k \cdot E[D_{\boldsymbol{\lambda}} \log p_k^*(X_k)] = \mathbf{0}.$$

Next, the components of the covariance matrix $\boldsymbol{\Sigma}(\boldsymbol{\lambda}) := \mathrm{Cov}(D_{\boldsymbol{\lambda}}\ell(\boldsymbol{\lambda}))$ are given by

$$(5.7) \qquad \Sigma_{st}(\boldsymbol{\lambda}) = \sum_k n_k \cdot \mathrm{Cov}(D_{\lambda_s} \log p_k^*(X_k), D_{\lambda_t} \log p_k^*(X_k))$$

and for the partial second derivatives we get

$$(5.8) \qquad J_{st}(\boldsymbol{\lambda}) := -D^2_{\lambda_s \lambda_t}\ell(\boldsymbol{\lambda}) = -\sum_k \sum_i D^2_{\lambda_s \lambda_t} \log p_k^*(X_{ki})$$

with expectation (cf. Appendix A.2)

$$(5.9) \quad I_{st}(\boldsymbol{\lambda}) := E(J_{st}(\boldsymbol{\lambda})) = \sum_k n_k \cdot E(D_{\lambda_s} \log p_k^*(X_k) \cdot D_{\lambda_t} \log p_k^*(X_k)).$$

Since $\ell(\boldsymbol{\lambda})$ is not the log-likelihood for sampling conditional on $X$, the matrices $\boldsymbol{\Sigma}(\boldsymbol{\lambda})$ and $\mathbf{I}(\boldsymbol{\lambda})$ need not be equal, but from (5.7) their difference is

$$(5.10) \quad I_{st}(\boldsymbol{\lambda}) - \Sigma_{st}(\boldsymbol{\lambda}) = \sum_k n_k \cdot E(D_{\lambda_s} \log p_k^*(X_k)) \cdot E(D_{\lambda_t} \log p_k^*(X_k)).$$

From now on we assume the essential:

CONDITION (R2). $\boldsymbol{\Sigma}(\boldsymbol{\lambda}) = \mathrm{Cov}(D_{\boldsymbol{\lambda}}\ell(\boldsymbol{\lambda}))$ is positive definite for all $\boldsymbol{\lambda}$.

Two equivalent formulations (cf. Appendix A.2) are

CONDITION (R2′). $\mathbf{I}(\boldsymbol{\lambda})$ is positive definite for all $\boldsymbol{\lambda}$.

CONDITION (R2″). For all $\boldsymbol{\theta}$, all $\mathbf{s} \in \mathbb{R}^S$ and $c_1, \ldots, c_K \in \mathbb{R}$: $D_{\boldsymbol{\theta}}\psi_{\boldsymbol{\theta}}(X, y_k) \cdot \mathbf{s} = c_k$ for $k = 1, \ldots, K$ almost surely $\Rightarrow \mathbf{s} = \mathbf{0}$.

In the last formulation—which does not include the nuisance parameter $\boldsymbol{\gamma}^*$—we can replace $X$ by $X_k$, since their distributions belong to $\mathcal{P}^X$ and hence dominate each other.

Using the block notation for an $(S+K) \times (S+K)$ matrix, say

$$\boldsymbol{\Sigma} = \begin{bmatrix} \boldsymbol{\Sigma}_{\boldsymbol{\theta\theta}} & \boldsymbol{\Sigma}_{\boldsymbol{\theta\gamma}} \\ \boldsymbol{\Sigma}_{\boldsymbol{\gamma\theta}} & \boldsymbol{\Sigma}_{\boldsymbol{\gamma\gamma}} \end{bmatrix},$$

a fundamental result can be derived (cf. Appendix A.2) by adopting the method in [12].

THEOREM 3. *For any $\boldsymbol{\lambda}$*

(a) $$\mathbf{I}(\boldsymbol{\lambda}) - \boldsymbol{\Sigma}(\boldsymbol{\lambda}) = \mathbf{I}(\boldsymbol{\lambda}) \cdot \begin{bmatrix} \mathbf{0} & \mathbf{0} \\ \mathbf{0} & \mathbf{W} \end{bmatrix} \cdot \mathbf{I}(\boldsymbol{\lambda}),$$

where the $K \times K$-matrix $\mathbf{W}$ is the sum of the diagonal $\mathrm{diag}(n_1^{-1}, \ldots, n_K^{-1})$ and the constant matrix $(n_0^{-1})$, that is, $W_{kl} = \Delta_{kl} n_k^{-1} + n_0^{-1}$ with the Kronecker's $\Delta$.

(b) $$[\mathbf{I}^{-1}(\boldsymbol{\lambda})]_{\boldsymbol{\theta\theta}} = [\mathbf{I}^{-1} \cdot \boldsymbol{\Sigma}(\boldsymbol{\lambda}) \cdot \mathbf{I}^{-1}(\boldsymbol{\lambda})]_{\boldsymbol{\theta\theta}}.$$

The matrix in (b) will later turn out to be the asymptotic covariance matrix of the estimate $\hat{\boldsymbol{\theta}}$.

*Log-bilinear association*: Using (4.5) and $\boldsymbol{\theta}$ (instead of $\tilde{\boldsymbol{\theta}}$) the model states

(5.11) $\quad \psi_{\boldsymbol{\theta}}(x, y_k) = \mathbf{z}^T \boldsymbol{\theta}_k \qquad$ with $\mathbf{z} = h_X(x),\ \boldsymbol{\theta} = (\boldsymbol{\theta}_1, \ldots, \boldsymbol{\theta}_K) \in \mathbb{R}^{K_X \times K}$

and is equivalent to the linear logistic regression model given by (5.2), that is,

$$\mathrm{logit}\, p_k^*(x) = \gamma_k^* + \mathbf{z}^T \boldsymbol{\theta}_k, \qquad k = 1, \ldots, K.$$

Condition (R2″) holds if $h_X(X)$ is not concentrated on a hyperplane of $\mathbb{R}^{K_X}$, that is, if the following condition is met (cf. Appendix A.2):

CONDITION $(\mathrm{R2})_{\mathrm{LBA}}$. For all $\mathbf{s} \in \mathbb{R}^{K_X} : \mathbf{s}^T h_X(X)$ is constant almost surely $\Rightarrow \mathbf{s} = \mathbf{0}$.

**6. Asymptotics and consistency.** We now turn to the asymptotic properties of the estimate $\hat{\boldsymbol{\lambda}} = (\hat{\boldsymbol{\theta}}, \hat{\boldsymbol{\gamma}}^*)$ in the model (3.13). Our asymptotic approach assumes that set $\Omega_Y^* = \{y_0, \ldots, y_K\}$ of conditional values will remain *fixed* while all subsample sizes $n_k$ tend to infinity with *fixed* ratios $\overline{r}_k = n_k/n > 0$ for all $n$ and $k$. Hence the nuisance parameter $\boldsymbol{\gamma}^*$, the distribution $P^*$ and its conditional densities $p_k^*(x)$ do not vary with $n$. The true parameter will now be denoted by $\boldsymbol{\lambda}^\circ = (\boldsymbol{\theta}^\circ, \boldsymbol{\gamma}^\circ)$ instead of $\boldsymbol{\lambda}$ and the notation $E, P$, etc. now refer to expectations, probabilities, etc. with respect to $\boldsymbol{\lambda}^\circ$. The conditional log-likelihood $\ell^{(n)}(\boldsymbol{\lambda})$—the additional index $n$ is supplied if necessary—need not have a unique maximizing argument $\hat{\boldsymbol{\lambda}}$ for *every* sample. Concerning uniqueness, the strong law of large numbers yields for the matrix $\mathbf{J}^{(n)}(\boldsymbol{\lambda}) = -D^2_{\boldsymbol{\lambda\lambda}} \ell^{(n)}(\boldsymbol{\lambda})$ from (5.8)

(6.1) $$\frac{1}{n} \mathbf{J}^{(n)}(\boldsymbol{\lambda}) \xrightarrow[n \to \infty]{} \overline{\mathbf{I}}(\boldsymbol{\lambda}) := \sum_{k=0}^{K} \overline{r}_k \cdot E(-D^2_{\boldsymbol{\lambda\lambda}} \log p_k^*(X_k)) \qquad \text{almost surely.}$$

The matrix $\overline{\mathbf{I}}(\boldsymbol{\lambda}) = \frac{1}{n}\mathbf{I}(\boldsymbol{\lambda})$ is positive definite by Condition (R2′) which implies $-D^2_{\boldsymbol{\lambda\lambda}} \ell^{(n)}(\boldsymbol{\lambda}) = -\mathbf{J}^{(n)}(\boldsymbol{\lambda})$ is negative definite for *almost all* (i.e., all except finitely many) $n$, almost surely. Hence—almost surely—the function $\ell^{(n)}(\boldsymbol{\lambda})$ is strictly concave for almost all $n$, which implies that $D_{\boldsymbol{\lambda}}\ell^{(n)}(\boldsymbol{\lambda}) = \mathbf{0}$ has at most one solution $\hat{\boldsymbol{\lambda}}$, which also maximizes $\ell^{(n)}(\boldsymbol{\lambda})$. Since the unique existence of a maximizing argument $\hat{\boldsymbol{\lambda}}$ of $\ell^{(n)}(\boldsymbol{\lambda})$ is not guaranteed for *every* $n$, we consider *any* sequence of (measurable) functions $\hat{\boldsymbol{\lambda}}^{(n)}$ as estimators if the estimating condition is met:

CONDITION (C1). If $\ell^{(n)}(\boldsymbol{\lambda})$ has a maximizing argument $\boldsymbol{\lambda}$, then $\ell^{(n)}(\hat{\boldsymbol{\lambda}}^{(n)}) = \operatorname{Max}_{\boldsymbol{\lambda}} \ell^{(n)}(\boldsymbol{\lambda})$.

To establish the consistency of such a sequence $\hat{\boldsymbol{\lambda}}^{(n)}$ we assume an integrability and an identifiability condition:

CONDITION (C2). $E\{\tilde{\psi}_X(X_k)\} < \infty$ for all $k = 0, \ldots, K$.

CONDITION (C3). $\psi_{\boldsymbol{\theta}_1}(X, y_k) = \psi_{\boldsymbol{\theta}_2}(X, y_k)$ for $k = 1, \ldots, K$ almost surely $\Rightarrow \boldsymbol{\theta}_1 = \boldsymbol{\theta}_2$.

As in Condition (R2″), we can equivalently replace $X$ by $X_k$ in Condition (C3). In Appendix A.3 we derive the asymptotic (unique) existence and strong consistency of the estimator:

THEOREM 4 (Consistency). *Under Conditions (C1)–(C3) the following properties hold almost surely:*

(a) *For almost all $n$ there exists a unique $\boldsymbol{\lambda}$ maximizing $\ell^{(n)}(\boldsymbol{\lambda})$, namely $\hat{\boldsymbol{\lambda}}^{(n)}$.*
(b) *For almost all $n$ there exists a unique solution $\boldsymbol{\lambda}$ of $D_{\boldsymbol{\lambda}}\ell^{(n)}(\boldsymbol{\lambda}) = \mathbf{0}$, namely $\hat{\boldsymbol{\lambda}}^{(n)}$.*
(c) $\hat{\boldsymbol{\lambda}}^{(n)} = (\hat{\boldsymbol{\theta}}^{(n)}, \boldsymbol{\gamma}^{*(n)}) \xrightarrow[n\to\infty]{} \boldsymbol{\lambda}^\circ = (\boldsymbol{\theta}^\circ, \boldsymbol{\gamma}^\circ).$

*Log-bilinear association*: In view of $\tilde{\psi}_X(x) = \|h_X(x)\|^2$, Condition (C2) reduces to a moment condition for $\mathbf{Z}_k = h_X(X_k)$:

CONDITION (C2)$_{\text{LBA}}$. $E\{\|\mathbf{Z}_k\|^2\} < \infty$ for all $k = 0, \ldots, K$.

And, using the parametrization (5.11), Condition (C3) reduces to

$$h_X^T(X)\boldsymbol{\theta}_{k1} = h_X^T(X)\boldsymbol{\theta}_{k2}$$

for $k = 1, \ldots, K$ almost surely $\quad \Rightarrow \quad \boldsymbol{\theta}_{1k} = \boldsymbol{\theta}_{2k}$ for all $k$,

which is implied by the stronger Condition (R2)$_{\text{LBA}}$.

**7. Asymptotic normality.** Let us finally establish the asymptotic normality for a sequence $\hat{\boldsymbol{\lambda}}^{(n)}$ of estimates. Instead of assuming Condition (C1), we derive the asymptotic distribution for *any* weakly consistent sequence $\hat{\boldsymbol{\lambda}}^{(n)}$ solving the estimating equation at least approximately, that is, we only assume

CONDITION (N1). $D_{\boldsymbol{\lambda}}\ell^{(n)}(\hat{\boldsymbol{\lambda}}^{(n)}) = o_P(\sqrt{n})$, respectively, $n^{-1/2} \cdot D_{\boldsymbol{\lambda}}\ell^{(n)} \times (\hat{\boldsymbol{\lambda}}^{(n)}) \xrightarrow[n \to \infty]{P} \mathbf{0}$.

CONDITION (N2). $\hat{\boldsymbol{\lambda}}^{(n)} \xrightarrow[n \to \infty]{P} \boldsymbol{\lambda}^\circ$.

Obviously both conditions hold under the assumptions of Theorem 4. Furthermore we assume the following consistency results, which are derived later (Theorem 6) from Condition (N2) and additional moment conditions:

CONDITION (N3). $\frac{1}{n}\int_0^1 \mathbf{J}^{(n)}(\boldsymbol{\lambda}^\circ + t[\hat{\boldsymbol{\lambda}}^{(n)} - \boldsymbol{\lambda}^\circ])\,dt \xrightarrow[n \to \infty]{P} \overline{\mathbf{I}}(\boldsymbol{\lambda}^\circ)$.

CONDITION (N4). $\frac{1}{n}\mathbf{J}^{(n)}(\hat{\boldsymbol{\lambda}}^{(n)}) \xrightarrow[n \to \infty]{P} \overline{\mathbf{I}}(\boldsymbol{\lambda}^\circ)$.

In Appendix A.4 we derive the asymptotic normality of the estimate as follows, where $\mathbf{A}^{1/2}$, respectively, $\mathbf{A}^{1/2}$ denotes the generalized Moore–Penrose inverse, respectively, the symmetric root of a positive semidefinite matrix $\mathbf{A}$, and $\mathbb{I}$ is the identity matrix.

THEOREM 5 (Normality). *Any sequence $\hat{\boldsymbol{\lambda}}^{(n)}$ of estimators with Conditions (N1)–(N3) is asymptotic normal*

(a) $\sqrt{n}[\hat{\boldsymbol{\lambda}}^{(n)} - \boldsymbol{\lambda}^\circ] \xrightarrow[n \to \infty]{\mathcal{L}} N(\mathbf{0}, \overline{\mathbf{I}}^{-1}(\boldsymbol{\lambda}^\circ) \cdot \overline{\boldsymbol{\Sigma}}(\boldsymbol{\lambda}^\circ) \cdot \overline{\mathbf{I}}^{-1}(\boldsymbol{\lambda}^\circ))$ *with* $\overline{\boldsymbol{\Sigma}}(\boldsymbol{\lambda}) := \sum_k \overline{r}_k \cdot \mathrm{Cov}(D_{\boldsymbol{\lambda}} \log p_k^*(X_k))$,

(b) $\sqrt{n}[\hat{\boldsymbol{\theta}}^{(n)} - \boldsymbol{\theta}^\circ] \xrightarrow[n \to \infty]{\mathcal{L}} N(\mathbf{0}, [\overline{\mathbf{I}}^{-1}(\boldsymbol{\lambda}^\circ)]_{\boldsymbol{\theta\theta}})$.

COROLLARY. *If in addition Condition (N4) holds, then*

(c) $([\mathbf{J}^{(n)}(\hat{\boldsymbol{\lambda}}^{(n)})^-]_{\boldsymbol{\theta\theta}}^{1/2})^- [\hat{\boldsymbol{\theta}}^{(n)} - \boldsymbol{\theta}^\circ] \xrightarrow[n \to \infty]{\mathcal{L}} N(\mathbf{0}, \mathbb{I})$.

Less formally (a) and (b) state

$$\hat{\boldsymbol{\lambda}} \underset{\text{as.}}{\sim} N(\boldsymbol{\lambda}^\circ, \mathbf{I}^{-1}(\boldsymbol{\lambda}^\circ) \cdot \boldsymbol{\Sigma}(\boldsymbol{\lambda}^\circ) \cdot \mathbf{I}^{-1}(\boldsymbol{\lambda}^\circ)), \qquad \hat{\boldsymbol{\theta}} \underset{\text{as.}}{\sim} N(\boldsymbol{\theta}, [\mathbf{I}^{-1}(\boldsymbol{\lambda}^\circ)]_{\boldsymbol{\theta\theta}}).$$

$\mathbf{J}(\hat{\boldsymbol{\lambda}})$ is a consistent estimate of $\mathbf{I}(\boldsymbol{\lambda}^\circ)$ by Condition (N4), and will be positive definite for almost all $n$ (almost surely) by (6.1). In this case, (c) states

$$\hat{\boldsymbol{\theta}} \underset{\text{as.}}{\sim} N(\boldsymbol{\theta}, [\mathbf{I}^{-1}(\hat{\boldsymbol{\lambda}})]_{\boldsymbol{\theta}\boldsymbol{\theta}}). \tag{7.1}$$

Notice that for an observed data set, the estimated covariance matrix $[\mathbf{J}^{-1}(\hat{\boldsymbol{\lambda}})]_{\boldsymbol{\theta}\boldsymbol{\theta}}$ (where the random variables are replaced by observations) is *identical* to the corresponding matrix under sampling conditional on $X$ (instead of $Y$). In this sense the estimate $\hat{\boldsymbol{\theta}}$ and its estimated asymptotic normal distribution are invariant under sampling conditional on either $Y$ or $X$. Hence asymptotic inference (i.e., tests or confidence regions) for the association parameter $\boldsymbol{\theta}$ based on the asymptotic distribution (7.1) of the estimate $\hat{\boldsymbol{\theta}}$ is invariant under both conditional sampling schemes, too.

The above Conditions (N3) and (N4) will now be derived from the consistency Condition (N2) and additional properties of the function $G$. For

$$H_r(\mathbf{z} \mid \boldsymbol{\theta}) = \sum_{k=0}^{K} |D_{\theta_r} G(\mathbf{z}, h_Y(y_k), \boldsymbol{\theta})|,$$

$$H_{rs}(\mathbf{z} \mid \boldsymbol{\theta}) = \sum_{k=0}^{K} |D^2_{\theta_r \theta_s} G(\mathbf{z}, h_Y(y_k), \boldsymbol{\theta})|,$$

$$H_{rst}(\mathbf{z} \mid \boldsymbol{\theta}) = \sum_{k=0}^{K} |D^3_{\theta_r \theta_s \theta_t} G(\mathbf{z}, h_Y(y_k), \boldsymbol{\theta})|,$$

the following result is proved in Appendix A.4.

THEOREM 6. *Conditions (N3) and (N4) follow from (N2) and the moment condition* $(\mathbf{MC})_{\mathrm{LBA}}$

CONDITION (MC). *There exists $\varepsilon^\circ > 0$ such that for $B(\boldsymbol{\theta}^\circ) = \{\boldsymbol{\theta} \mid \|\boldsymbol{\theta} - \boldsymbol{\theta}^\circ\| \leq \varepsilon^\circ\}$ and all $k = 0, \ldots, K$ the following functions of $\mathbf{Z}_k = h_X(X_k)$:*

$$\sup_{\boldsymbol{\theta} \in B(\boldsymbol{\theta}^\circ)} H_r(\mathbf{Z}_k \mid \boldsymbol{\theta})^3, \qquad \sup_{\boldsymbol{\theta} \in B(\boldsymbol{\theta}^\circ)} H_{st}(\mathbf{Z}_k \mid \boldsymbol{\theta})^2, \qquad \sup_{\boldsymbol{\theta} \in B(\boldsymbol{\theta}^\circ)} H_{rst}(\mathbf{Z}_k \mid \boldsymbol{\theta})$$

*have finite expectation for all $r, s, t = 1, \ldots, S$.*

Hence the requirements for Theorem 5 are met if Conditions (MC) and (C1)–(C3) in Theorem 4 hold.

*Log-bilinear association*: The log-bilinear association model is based on the function $G(\mathbf{z}, \mathbf{v}, \boldsymbol{\theta}) = \mathbf{z}^T \boldsymbol{\theta} \mathbf{v}$ with partial derivatives $D_{\theta_{lm}} G(\mathbf{z}, \mathbf{v}, \boldsymbol{\theta}) = z_l v_m$ and vanishing higher derivatives. Hence Condition (MC) holds if Condition $(\mathrm{C2})_{\mathrm{LBA}}$ is strengthened to

CONDITION $(\mathrm{MC})_{\mathrm{LBA}}$. $E\{\|\mathbf{Z}_k\|^3\} < \infty$ for all $k = 0, \ldots, K$.

**8. Discussion.** Association models for a pair of random elements $(X, Y)$ do not restrict the marginal distributions of $X$ and $Y$ but only their odds ratio function. We have looked at parametric association models which include the important log-bilinear association models. An advantage of these models is that inference about the odds ratio (or association) parameter vector $\boldsymbol{\theta}$ may be obtained from sampling $Y$ conditional on fixed values of $X$ or vice versa. The maximum likelihood estimate $\hat{\boldsymbol{\theta}}$ is the same under both conditional sampling schemes, and asymptotic inference concerning $\boldsymbol{\theta}$ is invariant with respect to sampling, too. More precisely, we have shown that for samples conditional on $Y$, the estimate $\hat{\boldsymbol{\theta}}$ maximizing the "reverse" conditional likelihood $L_{Y|X}$ is consistent, asymptotic normal and its estimated asymptotic covariance matrix is the same as if sampling had been conditional on $X$. These results have been obtained much earlier for *discrete* $Y$ with *finite* range for the multivariate linear logistic regression model in [12] and for the general logistic regression model in [16] (for $X$ with finite range) and [14]. Our result allows both $X$ and $Y$ to be arbitrary random vectors each having discrete and/or continuous components.

Furthermore, asymptotic inference for the regression parameters $\boldsymbol{\beta}$ in widely used regression models is available when sampling is conditional on $Y$ (instead of $X$). For example, in log-linear regression models for Poisson variates we have $\boldsymbol{\beta} = \boldsymbol{\theta}$ and hence inference on $\boldsymbol{\beta}$ may also be obtained from samples conditional on $Y$. Even in the linear regression model $\mu(x) = \alpha + \mathbf{z}^T\boldsymbol{\beta}$ with covariate vector $\mathbf{z} = h_X(x)$ and $\mathcal{L}(Y \mid x) = N(\mu(x), \sigma^2)$, asymptotic inference for $\boldsymbol{\theta} = \sigma^{-2}\boldsymbol{\beta}$ may be obtained from samples conditional on $Y$—including tests of a linear hypothesis $\mathbf{C}\boldsymbol{\theta} = \mathbf{0}$, which is equivalent to $\mathbf{C}\boldsymbol{\beta} = \mathbf{0}$. However, confidence regions are only available for $\boldsymbol{\theta}$, but not for $\boldsymbol{\beta}$, unless an estimate of $\sigma^2$ from another sample is at hand. This extends to the multivariate case where the conditional distribution of $Y$ is multivariate normal $N_K(\mu(x), \boldsymbol{\Sigma})$ and the odds ratio parameter is given by $\boldsymbol{\theta} = \boldsymbol{\beta}\boldsymbol{\Sigma}^{-1}$. Although sampling conditional on $Y$ seems unnatural for a regression model, it may be very attractive if such a sample is much easier (e.g., cheaper or quicker) to obtain. The advantages of (retrospective) case-control over (prospective) cohort studies can thus be extended to an arbitrary response vector $Y$, for example, to *infinite discrete* response categories or to a *continuous* response $Y$. In the latter case we do not get confidence intervals for $\boldsymbol{\beta}$, but *tests* for linear hypothesis—which may be of primary interest (e.g., in a clinical trial)—are available.

Related, but different, semiparametric models for random vectors $X = (X_1, \ldots, X_I)$ and $Y = (Y_1, \ldots, Y_J)$ are given by multivariate copulas which specify parametric distributions on $[0, 1]^{I+J}$ with uniform marginals. However, a copula is not an association model in our sense (cf. [9]) because a copula only leaves the marginal distributions of all *univariate components* $X_i$ and $Y_j$ arbitrary, but the marginal distribution of the *vectors* $X$, respectively, $Y$ are restricted through the parametrization of the copula, unless both $X$ and $Y$ are univariate. And even in the latter case, the odds ratio function $OR(X, Y)$ cannot be recovered from the corresponding copula unless both marginal distributions of $X$ and $Y$ are *known*. Hence the rather general semiparametnc associations models considered here do not fit in the framework of copulas.

## APPENDIX: PROOFS

**A.1. Proof of Theorem 2 (existence).** We have already seen that (b) implies (c) and it remains to derive (a) from (c), which uses the concept of an $I$-projection and heavily relies on results by Csiszár [2] and Rüschendorf and Thomsen [13]. Setting $\pi = \pi_X \times \pi_Y$ we first conclude from Condition (E2) the existence of $R \in \mathcal{P}$ with $\pi$-density

$$r = \exp(\psi - \beta - \gamma - \alpha) > 0, \qquad \alpha = \log \int \exp(\psi - \beta - \gamma) \, d\pi$$

and the wanted $P$ will be the $I$-projection of $R$ on $\mathcal{E} = \{P \in \mathcal{P} \mid P^X = \pi_X, P^Y = \pi_Y\}$. The integrability of $\psi$, $\beta$ and $\gamma$ implies

$$I(\pi \mid R) = \int \log\left(\frac{1}{r}\right) d\pi = \int (\alpha + \beta + \gamma - \psi) \, d\pi < \infty$$

and since $\pi \in \mathcal{E}$, we conclude from Theorem 2.1 in [2] that $R$ has an $I$-projection $P$ on $\mathcal{E}$. Application of Theorem 3.1 in [2] to the set

$$\mathcal{F} = \{f_X + f_Y \mid f_X \in \mathcal{L}_1(\pi_X), f_Y \in \mathcal{L}_1(\pi_Y)\} \subset \mathcal{L}_1(P)$$

yields that the $R$-density $p_R$ of $P$ satisfies $p_R = \exp(h)$ $\pi$-almost surely, where $h$ belongs to the closure $\mathcal{F}^-$ of $\mathcal{F}$ in $\mathcal{L}_1(P)$. Rüschendorf and Thomsen [13] pointed out that $\mathcal{F}$ need not be closed in $\mathcal{L}_1(P)$—which was claimed in the proof of Corollary 3.1, case (B) in Csiszár [2].

Now $R \ll \pi$ implies that $\exp(h) > 0$ is an $R$-density of $P$ and hence $R \ll P \ll R$. Furthermore $r > 0$ yields $R \ll \pi \ll R$ and hence $P \in \mathcal{P}_{\ll}$, since $P^{XY} = \pi$. From Theorem 2.2 in [2] we obtain

$$I(\pi \mid P) + I(P \mid R) \leq I(\pi \mid R) < \infty,$$

which establishes $P \in \mathcal{P}_f$. Finally $OR(P) = \psi$ remains to be shown. From $P \ll P^{XY}$ and Proposition 2 in [13] we conclude the existence of measurable functions $a: \Omega_X \to \mathbb{R}$ and $b: \Omega_Y \to \mathbb{R}$, such that $h(x, y) = b(x) + c(y)$ $P$-almost surely, and hence $R$-almost surely. Hence a $\pi$-density of $p$ is given by

$$\frac{dP}{d\pi} = \frac{dP}{dR} \cdot \frac{dR}{d\pi} = \exp(b + c) \cdot r = \exp(b + c - \beta - \gamma - \alpha + \psi)$$

and a direct calculation yields $\log OR(P) = \psi$ as required.

**A.2. Proof of the results in Section 5.** We start with some preliminary results. The derivatives of $\log p_k^*$ are given by

(A.1)
$$D_{\lambda_s} \log p_k^*(x) = \frac{D_{\lambda_s} p_k^*(x)}{p_k^*(x)},$$
$$D_{\lambda_s \lambda_t}^2 \log p_k^*(x) = \frac{D_{\lambda_s \lambda_t}^2 p_k^*(x)}{p_k^*(x)} - D_{\lambda_s} \log p_k^*(x) \cdot D_{\lambda_t} \log p_k^*(x).$$

For any set of measurable functions $G_k(x)$ we obtain from (5.1) a key equality:

(A.2)
$$\sum_k \overline{r}_k \cdot E(G_k(X_k)) = \sum_k \overline{r}_k \cdot E(G_k(X) \mid Y = y_k)$$
$$= \sum_k \overline{r}_k \cdot \int G_k(x) \cdot p(x \mid Y = y_k) \, d\nu_X(x)$$
$$= \int \sum_k G_k(x) \cdot p_k^*(x) \cdot p^{*X}(x) \, d\nu_X(x)$$
$$= E^*\left[\sum_k G_k(X) \cdot p_k^*(X)\right],$$

where $E^*$ denotes expectation with respect to $P^*$.

In particular, we get for $G_k(x) = H(x) \cdot D_{\boldsymbol{\lambda}} \log p_k^*(x)$ and any measurable $H(x)$

(A.3)
$$\sum_k \overline{r}_k \cdot E[H(X_k) \cdot D_{\boldsymbol{\lambda}} \log p_k^*(X_k)]$$
$$= E^*\left[\sum_k H(X) \cdot D_{\boldsymbol{\lambda}} \log p_k^*(X) \cdot p_k^*(X)\right]$$
$$= E^*\left[H(X) \cdot \sum_k D_{\boldsymbol{\lambda}} p_k^*(X)\right] = \mathbf{0},$$

since $p_+^*(x) = 1$. In particular, (5.6) follows for $H(x) = 1$.

PROOF OF (5.9). Choosing $G_k(X_k) = p_k^*(X_k)^{-1} \cdot D_{\lambda_s \lambda_t}^2 p_k^*(X_k)$ in (A.2) yields

$$\sum_k \overline{r}_k \cdot E[p_k^*(X_k)^{-1} \cdot D_{\lambda_s \lambda_t}^2 p_k^*(X_k)] = E^*\left[\sum_k D_{\lambda_s \lambda_t}^2 p_k^*(X_k)\right] = \mathbf{0}$$

and (5.9) follows using (A.1):

$$E(J_{st}(\boldsymbol{\lambda})) = n \cdot \sum_k \overline{r}_k \cdot E(D_{\lambda_s} \log p_k^*(X_k) \cdot D_{\lambda_t} \log p_k^*(X_k)).$$

$\square$

PROOF OF CONDITIONS (R2) $\Leftrightarrow$ (R2'). By (5.10) $\mathbf{I}(\boldsymbol{\lambda})$ is a sum of $\boldsymbol{\Sigma}(\boldsymbol{\lambda})$ and a positive semidefinite matrix. Hence $\mathbf{I}(\boldsymbol{\lambda})$ is positive semidefinite, and even positive definite, provided Condition (R2) holds. Conversely, let Condition (R2') hold. Then $\mathbf{t}^T \boldsymbol{\Sigma}(\boldsymbol{\lambda}) \mathbf{t} = \mathrm{Var}(\mathbf{t}^T D_{\boldsymbol{\lambda}} \ell(\boldsymbol{\lambda})^T) = 0$ implies that $\mathbf{t}^T D_{\boldsymbol{\lambda}} \ell(\boldsymbol{\lambda})^T$ is constant almost surely, and hence $\mathbf{t}^T D^2_{\boldsymbol{\lambda}\boldsymbol{\lambda}} \ell(\boldsymbol{\lambda}) = D_{\boldsymbol{\lambda}}[\mathbf{t}^T D_{\boldsymbol{\lambda}} \ell(\boldsymbol{\lambda})^T] = \mathbf{0}$ almost surely. Thus $\mathbf{t}^T \mathbf{I}(\boldsymbol{\lambda}) = E(\mathbf{t}^T D^2_{\boldsymbol{\lambda}\boldsymbol{\lambda}} \ell(\boldsymbol{\lambda})) = \mathbf{0}$, which implies $\mathbf{t} = \mathbf{0}$ by Condition (R2'). Hence Condition (R2) holds. $\square$

PROOF OF CONDITIONS (R2') $\Leftrightarrow$ (R2''). $\mathbf{I}(\boldsymbol{\lambda})$ is positive semidefinite (as already observed) and hence Condition (R2') is equivalent to

(A.4) $\qquad$ For all $\mathbf{t} \in \mathbb{R}^{S+K}$ $\qquad \mathbf{t}^T \mathbf{I}(\boldsymbol{\lambda}) \mathbf{t} = 0 \quad \Rightarrow \quad \mathbf{t} = \mathbf{0}$.

For any $\mathbf{t} \in \mathbb{R}^{S+K}$ we get from (5.9)

$$\mathbf{t}^T \mathbf{I}(\boldsymbol{\lambda}) \mathbf{t} = \sum_k n_k \cdot E(\|D_{\boldsymbol{\lambda}} \log p_k^*(X_k) \cdot \mathbf{t}\|^2)$$

and since the distributions of $X_k$ and $X$ dominate each other:

(A.5) $\qquad \mathbf{t}^T \mathbf{I}(\boldsymbol{\lambda}) \mathbf{t} = 0 \quad \Leftrightarrow \quad D_{\boldsymbol{\lambda}} \log p_k^*(X) \cdot \mathbf{t} = 0$

$\qquad\qquad\qquad\qquad\qquad\qquad\qquad$ for $k = 0, \ldots, K$ almost surely.

To derive Condition (R2') from Condition (R2''), let $\mathbf{t}^T \mathbf{I}(\boldsymbol{\lambda}) \mathbf{t} = 0$. From (5.4) we get

(A.6) $\qquad \mathrm{logit}\, p_k^*(X) = \log p_k^*(X) - \log p_0^*(X) = \gamma_k^* + \psi_{\boldsymbol{\theta}}(X, y_k)$

and for $\mathbf{t} = (\mathbf{s}, -\mathbf{c})$ with $\mathbf{s} \in \mathbb{R}^S$, $\mathbf{c} = (c_1, \ldots, c_K)$, we obtain from (A.5) almost surely

(A.7) $\begin{aligned} 0 &= D_{\boldsymbol{\lambda}} \mathrm{logit}\, p_k^*(X) \cdot \mathbf{t} = D_{\boldsymbol{\theta}} \mathrm{logit}\, p_k^*(X) \cdot \mathbf{s} - D_{\boldsymbol{\gamma}} \mathrm{logit}\, p_k^*(X) \cdot \mathbf{c} \\ &= D_{\boldsymbol{\theta}} \psi_{\boldsymbol{\theta}}(X, y_k) \cdot \mathbf{s} - c_k \qquad \text{for all } k = 1, \ldots, K. \end{aligned}$

And from Condition (R2'') we conclude $\mathbf{s} = \mathbf{0}$ as well as $c_k = 0$ for all $k$, and thus $\mathbf{t} = \mathbf{0}$.

Conversely, suppose Condition (R2') holds. To establish Condition (R2''), it suffices to show that (A.7) implies $\mathbf{s} = \mathbf{0}$. From (5.2) and (5.4) we get

$$p_0^*(X) = \left( \sum_l \exp[\mathrm{logit}\, p_l^*(X, y_l)] \right)^{-1},$$

$$D_{\boldsymbol{\lambda}} \log p_0^*(X) \cdot \mathbf{t} = p_0^*(X)^{-1} \sum_l \exp[\mathrm{logit}\, p_l^*(X, y_l)] \cdot D_{\boldsymbol{\lambda}} \mathrm{logit}\, p_l^*(X, y_l) \cdot \mathbf{t}.$$

Hence (A.7)—and $\operatorname{logit} p_0^* = 0$—imply $D_{\boldsymbol{\lambda}} \log p_0^*(X) \cdot \mathbf{t} = 0$ almost surely. From (A.6) we get $D_{\boldsymbol{\lambda}} \log p_k^*(X) \cdot \mathbf{t} = 0$ for $k = 0, \ldots, K$ almost surely, and (A.5), (A.4) establish $\mathbf{t} = \mathbf{0}$ and hence $\mathbf{s} = \mathbf{0}$. □

PROOF OF THEOREM 3. Part (a) is equivalent to three equations:

(a)$_{\boldsymbol{\theta\theta}}$ $\qquad \mathbf{I}_{\boldsymbol{\theta\theta}} - \boldsymbol{\Sigma}_{\boldsymbol{\theta\theta}} = \mathbf{I}_{\boldsymbol{\theta\gamma}} \cdot \mathbf{W} \cdot \mathbf{I}_{\boldsymbol{\theta\gamma}}^T,$

(a)$_{\boldsymbol{\theta\gamma}}$ $\qquad \mathbf{I}_{\boldsymbol{\theta\gamma}} - \boldsymbol{\Sigma}_{\boldsymbol{\theta\gamma}} = \mathbf{I}_{\boldsymbol{\theta\gamma}} \cdot \mathbf{W} \cdot \mathbf{I}_{\boldsymbol{\gamma\gamma}},$

(a)$_{\boldsymbol{\gamma\gamma}}$ $\qquad \mathbf{I}_{\boldsymbol{\gamma\gamma}} - \boldsymbol{\Sigma}_{\boldsymbol{\gamma\gamma}} = \mathbf{I}_{\boldsymbol{\gamma\gamma}} \cdot \mathbf{W} \cdot \mathbf{I}_{\boldsymbol{\gamma\gamma}}.$

Some prerequisite results are derived first using the notation

$$b_{sk} = E[D_{\theta_s} \log p_k^*(X_k)] \in \mathbb{R}, \qquad \mathbf{b}_k = (b_{1k}, \ldots, b_{Sk}) \in \mathbb{R}^S,$$

$$c_{mk} = E[D_{\gamma_m^*} \log p_k^*(X_k)] \in \mathbb{R}, \qquad \mathbf{c}_k = (c_{1k}, \ldots, c_{Kk}) \in \mathbb{R}^K,$$

$$\mathbf{B} = (\mathbf{b}_1, \ldots, \mathbf{b}_K) \in \mathbb{R}^{S \times K}, \qquad \overline{\mathbf{B}} = (\mathbf{b}_0, \ldots, \mathbf{b}_K) \in \mathbb{R}^{S \times (K+1)},$$

$$\mathbf{C} = (\mathbf{c}_1, \ldots, \mathbf{c}_K) \in \mathbb{R}^{K \times K}, \qquad \overline{\mathbf{C}} = (\mathbf{c}_0, \ldots, \mathbf{c}_K) \in \mathbb{R}^{K \times (K+1)},$$

$$\mathbf{N} = \operatorname{diag}(n_1, \ldots, n_K) \in \mathbb{R}^{K \times K}, \qquad \overline{\mathbf{N}} = \operatorname{diag}(n_0, \ldots, n_K) \in \mathbb{R}^{K \times (K+1)}.$$

From (5.3) we obtain the partial derivatives

$$D_{\theta_s} \log p_k^*(x) = D_{\theta_s} \psi_{\boldsymbol{\theta}}(x, y_k) - \sum_l p_l^*(x) \cdot D_{\theta_s} \psi_{\boldsymbol{\theta}}(x, y_l),$$

$$D_{\gamma_m^*} \log p_k^*(x) = \Delta_{km} - p_m^*(x)$$

and (5.9) yields

$$I_{\lambda_t \gamma_m^*} = n \sum_k \overline{r}_k \cdot E(D_{\theta_s} \log p_k^*(X_k) \cdot D_{\gamma_m^*} \log p_k^*(X_k))$$

$$= n_m \cdot E(D_{\theta_s} \log p_m^*(X_k)) - n \sum_k \overline{r}_k \cdot E(p_m^*(X_k) \cdot D_{\theta_s} \log p_k^*(X_k))$$

$$= n_m \cdot E(D_{\theta_s} \log p_m^*(X_k)) \qquad [\text{cf. (A.3) for } H(x) = p_m^*(x)].$$

Hence $I_{\theta_s \gamma_m^*} = n_m \cdot b_{sm}$, $I_{\gamma_l^* \gamma_m^*} = n_m \cdot c_{lm}$, or in matrix notation

(A.8) $\qquad \mathbf{I}_{\boldsymbol{\theta\gamma}} = \mathbf{B} \cdot \mathbf{N}, \qquad \mathbf{I}_{\boldsymbol{\gamma\gamma}} = \mathbf{C} \cdot \mathbf{N}.$

From (5.6) we have $\sum_k n_k \cdot b_{sk} = 0$ and $\sum_k n_k \cdot c_{mk} = 0$, or in matrix notation

(A.9) $\qquad \mathbf{0} = n_0 \mathbf{b}_0 + \mathbf{B}\mathbf{n}, \qquad \mathbf{0} = n_0 \mathbf{c}_0 + \mathbf{C}\mathbf{n}, \qquad \mathbf{n} = (n_1, \ldots, n_K).$

Using the constant vector $\mathbf{e}_+ = (1)$ and constant matrix $\mathbf{e}_+ \mathbf{e}_+^T = (1)$ we thus obtain

$$\mathbf{I}_{\boldsymbol{\theta\gamma}} \cdot \mathbf{W} = \mathbf{B} \cdot \mathbf{N}[n_0^{-1} \mathbf{e}_+ \mathbf{e}_+^T + \mathbf{N}^{-1}] = n_0^{-1} \mathbf{B} \cdot \mathbf{n} \cdot \mathbf{e}_+^T + \mathbf{B} = -\mathbf{b}_0 \cdot \mathbf{e}_+^T + \mathbf{B}$$

and similarly with $\mathbf{C}$ instead of $\mathbf{B}$
$$\mathbf{I}_{\gamma\gamma} \cdot \mathbf{W} = \mathbf{C} \cdot \mathbf{N} \cdot \mathbf{W} = -\mathbf{c}_0 \cdot \mathbf{e}_+^T + \mathbf{C}.$$

Now $(a)_{\boldsymbol{\theta}\gamma}$ is obtained as follows:

$$\begin{aligned}
\mathbf{I}_{\boldsymbol{\theta}\gamma} \cdot \mathbf{W} \cdot \mathbf{I}_{\gamma\gamma}^T &= [\mathbf{B} - \mathbf{b}_0 \cdot \mathbf{e}_+^T][\mathbf{C} \cdot \mathbf{N}]^T = \mathbf{B} \cdot \mathbf{N} \cdot \mathbf{C}^T - \mathbf{b}_0 \cdot [\mathbf{C} \cdot \mathbf{n}]^T & [\text{cf. (A.8)}] \\
&= \mathbf{B} \cdot \mathbf{N} \cdot \mathbf{C}^T + \mathbf{b}_0 \cdot n_0 \cdot \mathbf{c}_0^T = \overline{\mathbf{B}} \cdot \overline{\mathbf{N}} \cdot \overline{\mathbf{C}}^T & [\text{cf. (A.9)}] \\
&= \mathbf{I}_{\boldsymbol{\theta}\gamma} - \boldsymbol{\Sigma}_{\boldsymbol{\theta}\gamma} & [\text{cf. (5.10)}]
\end{aligned}$$

And $(a)_{\boldsymbol{\theta\theta}}$, respectively $(a)_{\gamma\gamma}$, is established similarly (replace $\mathbf{B}$ and $\mathbf{b}_0$ by $\mathbf{C}$ and $\mathbf{c}_0$, respectively, vice versa). Hence (a) holds, and multiplication with $\mathbf{I}^{-1}(\boldsymbol{\lambda})$ yields (b). □

PROOF OF CONDITION $(\text{R2})_{\text{LBA}} \Rightarrow (\text{R2}'')$. Suppose for $\mathbf{s} = (\mathbf{s}_1, \ldots, \mathbf{s}_K) \in \mathbb{R}^{K_X \times K}$ and $c_1, \ldots, c_K \in \mathbb{R}$ we have for all $k = 1, \ldots, K$

$$c_k = D_{\boldsymbol{\theta}} \psi_{\boldsymbol{\theta}}(X, y_k) \cdot \mathbf{s} = \sum_l D_{\boldsymbol{\theta}_l} \psi_{\boldsymbol{\theta}}(X, y_k) \cdot \mathbf{s}_l = h_X(X)^T \cdot \mathbf{s}_k \qquad \text{almost surely.}$$

Then Condition $(\text{R2})_{\text{LBA}}$ implies $\mathbf{s}_k = \mathbf{0}$ for all $k$, and hence $\mathbf{s} = \mathbf{0}$. □

**A.3. Proof of Theorem 4 (consistency).** The proof is based on the ingenious ideas from Wald [15]. The log-odds ratio $\psi_{\boldsymbol{\theta}}(x,y)$ in the model (3.13) depends only on the vectors $\mathbf{z} = h_X(x)$ and $\mathbf{v} = h_Y(y)$. Therefore we regard $p_k^*(x) = \tilde{p}_k(\mathbf{z} \mid \boldsymbol{\lambda})$ as a function of $\mathbf{z}$ and $\boldsymbol{\lambda}$ using the notation

$$G_k(\mathbf{z}, \boldsymbol{\theta}) := G(\mathbf{z}, h_Y(y_k), \boldsymbol{\theta}) = \psi_{\boldsymbol{\theta}}(x, y_k),$$

$$\tilde{p}_k(\mathbf{z} \mid \boldsymbol{\lambda}) := \frac{\exp[\gamma_k^* + G_k(\mathbf{z}, \boldsymbol{\theta})]}{\sum_l \exp[\gamma_l^* + G_l(\mathbf{z}, \boldsymbol{\theta})]} = p_k^*(x),$$

$$\eta_k(\mathbf{z} \mid \boldsymbol{\lambda}) := \log \tilde{p}_k(\mathbf{z} \mid \boldsymbol{\lambda}) = \gamma_k^* + G_k(\mathbf{z}, \boldsymbol{\theta}) - \log\left(\sum_l \exp[\gamma_l^* + G_l(\mathbf{z}, \boldsymbol{\theta})]\right).$$

We first show for $\mathbf{Z}_k := h_X(X_k)$

(A.10) $\qquad E\{|\eta_k(\mathbf{Z}_k \mid \boldsymbol{\lambda})|\} < \infty \qquad$ for all $\boldsymbol{\lambda}$ and $k = 0, \ldots, K$.

From $\gamma_0^* = 0 = G_0(\mathbf{z}, \boldsymbol{\theta})$ and $\tilde{p}_0(\mathbf{z} \mid \boldsymbol{\lambda}) \leq 1$ we get

$$|\eta_0(\mathbf{z} \mid \boldsymbol{\lambda})| = \log\left(\sum_l \exp[\gamma_l^* + G_l(\mathbf{z}, \boldsymbol{\theta})]\right) \leq \log(K+1) + \|\boldsymbol{\gamma}^*\| + \operatorname*{Max}_l |G_l(\mathbf{z}, \boldsymbol{\theta})|.$$

And Condition (OR2) yields

(A.11) $\qquad |G_l(\mathbf{z}, \boldsymbol{\theta})| \leq [\tilde{\psi}_X(x) + \tilde{\psi}_Y(y_l)] \cdot \|\boldsymbol{\theta}\|,$

which in view of Condition (C2) proves (A.10) for $k=0$. For $k>0$ we get

$$|\eta_k(\mathbf{z}\mid\boldsymbol{\lambda})| = |\gamma_k^* + G_k(\mathbf{z},\boldsymbol{\theta}) + \eta_0(\mathbf{z},\boldsymbol{\lambda})| \leq \|\boldsymbol{\gamma}^*\| + |G_k(\mathbf{z},\boldsymbol{\theta})| + |\eta_0(\mathbf{z}\mid\boldsymbol{\lambda})|.$$

Hence (A.11) and Condition (C2) establish (A.10).

Next we prove three basic lemmas.

LEMMA A.1. *For any* $\boldsymbol{\lambda} \neq \boldsymbol{\lambda}^\circ : \sum_{k=0}^{K} \overline{r}_k \cdot E\{\eta_k(\mathbf{Z}_k \mid \boldsymbol{\lambda}) - \eta_k(\mathbf{Z}_k \mid \boldsymbol{\lambda}^\circ)\} < 0.$

LEMMA A.2. *For* $k=0,\ldots,K$ *and any* $\boldsymbol{\lambda}$:

$$\lim_{\varepsilon \to 0} E\left\{\sup_{\|\boldsymbol{\lambda}'-\boldsymbol{\lambda}\|\leq\varepsilon} \eta_k(\mathbf{Z}_k \mid \boldsymbol{\lambda}')\right\} = E\{\eta_k(\mathbf{Z}_k \mid \boldsymbol{\lambda})\}.$$

LEMMA A.3. *For any compact set* $A \subset \mathbb{R}^K \times \mathbb{R}^S$ *with* $\boldsymbol{\lambda}^\circ \notin A$:

$$\lim_{n\to\infty}\left[\sup_{\boldsymbol{\lambda}\in A} \ell^{(n)}(\boldsymbol{\lambda}) - \ell^{(n)}(\boldsymbol{\lambda}^\circ)\right] = -\infty \qquad \text{almost surely.}$$

PROOF OF LEMMA A.1. $U_k = \eta_k(\mathbf{Z}_k \mid \boldsymbol{\lambda}) - \eta_k(\mathbf{Z}_k \mid \boldsymbol{\lambda}^\circ)$ has finite expectation by (A.10), and Jensen's inequality yields

$$\sum_k \overline{r}_k \cdot E\{U_k\} \leq \sum_k \overline{r}_k \cdot \log E\{\exp(U_k)\} \leq \log\left(\sum_k \overline{r}_k \cdot E\{\exp(U_k)\}\right).$$
(A.12)

Equation (A.2) with $G_k(X_k) = \exp(U_k) = \tilde{p}_k(\mathbf{Z}_k \mid \boldsymbol{\lambda})[\tilde{p}_k(\mathbf{Z}_k \mid \boldsymbol{\lambda}^\circ)]^{-1}$ and $\boldsymbol{\lambda} \neq \boldsymbol{\lambda}^\circ$ gives (the true parameter is denoted by $\boldsymbol{\lambda}^\circ$ here)

$$\sum_k \overline{r}_k \cdot E\{\exp(U_k)\} = E^*\left\{\sum_k \tilde{p}_k(\mathbf{Z}_k \mid \boldsymbol{\lambda})\right\} = 1$$

and (A.12) implies $\sum_k \overline{r}_k \cdot E\{U_k\} \leq 0$. It remains to show that this inequality is strict. Suppose not; then equality holds in both places of (A.12). The first equality implies that each $U_k$ is constant almost surely, say $U_k = \log c_k$, and the second yields $c_k = c$ for all $k$, hence $U_k = \log c$, respectively, $\tilde{p}_k(\mathbf{Z}_k \mid \boldsymbol{\lambda}) = c \cdot \tilde{p}_k(\mathbf{Z}_k \mid \boldsymbol{\lambda}^\circ)$ almost surely. From $\sum_k \tilde{p}_k = 1$ we get $c = 1$, and hence

(A.13) $\qquad \eta_k(\mathbf{Z}_k \mid \boldsymbol{\lambda}) = \eta_k(\mathbf{Z}_k \mid \boldsymbol{\lambda}^\circ) \qquad$ for all $k$ almost surely.

Then

$$\psi_{\boldsymbol{\theta}}(X_k, y_k) = \eta_k(\mathbf{Z}_k \mid \boldsymbol{\lambda}) + \eta_0(\mathbf{Z}_0 \mid \boldsymbol{\lambda}) - \eta_0(\mathbf{Z}_k \mid \boldsymbol{\lambda}) - \eta_k(\mathbf{Z}_0 \mid \boldsymbol{\lambda}) = \psi_{\boldsymbol{\theta}^\circ}(X_k, y_k)$$

almost surely, and since the distributions of $X_k$ and $X$ dominate each other,

$$\psi_{\boldsymbol{\theta}}(X, y_k) = \psi_{\boldsymbol{\theta}^\circ}(X, y_k) \qquad \text{for all } k \text{ almost surely.}$$

From Condition (C3) we get $\boldsymbol{\theta} = \boldsymbol{\theta}^\circ$. For $\boldsymbol{\lambda} = (\boldsymbol{\theta}, \boldsymbol{\gamma}^*)$ (A.13) gives almost surely

$$\gamma_k^* + G_k(\mathbf{Z}_k, \boldsymbol{\theta}) = \eta_k(\mathbf{Z}_k \mid \boldsymbol{\lambda}) - \eta_0(\mathbf{Z}_k \mid \boldsymbol{\lambda}) = \gamma_k^\circ + G_k(\mathbf{Z}_k, \boldsymbol{\theta}^\circ) \qquad \text{for all } k$$

and from $\boldsymbol{\theta} = \boldsymbol{\theta}^\circ$ we conclude $\boldsymbol{\gamma}^* = \boldsymbol{\gamma}^\circ$, which contradicts $\boldsymbol{\lambda} \neq \boldsymbol{\lambda}^\circ$.  □

PROOF OF LEMMA A.2. Continuity implies for any positive sequence $\varepsilon_n \to 0$

$$\sup_{\|\boldsymbol{\lambda}' - \boldsymbol{\lambda}\| \leq \varepsilon_n} \eta_k(\mathbf{z} \mid \boldsymbol{\lambda}') \xrightarrow[n\to\infty]{} \eta_k(\mathbf{z} \mid \boldsymbol{\lambda}).$$

Since

(A.14) $$\eta_k(\mathbf{z} \mid \boldsymbol{\lambda}) \leq \sup_{\|\boldsymbol{\lambda}' - \boldsymbol{\lambda}\| \leq \varepsilon_n} \eta_k(\mathbf{z} \mid \boldsymbol{\lambda}') \leq 0,$$

the dominated convergence theorem and (A.10) yield

$$E\left\{ \sup_{\|\boldsymbol{\lambda}' - \boldsymbol{\lambda}\| \leq \varepsilon_n} \eta_k(\mathbf{Z}_k \mid \boldsymbol{\lambda}') \right\} \xrightarrow[n\to\infty]{} E\{\eta_k(\mathbf{Z}_k \mid \boldsymbol{\lambda})\}. \qquad \square$$

PROOF OF LEMMA A.3. For $\varepsilon > 0$ consider the ball $B(\boldsymbol{\lambda} \mid \varepsilon) = \{\boldsymbol{\lambda}' \mid \|\boldsymbol{\lambda}' - \boldsymbol{\lambda}\| \leq \varepsilon\}$ with interior $B^\circ(\boldsymbol{\lambda} \mid \varepsilon)$ and let $\eta_k(\mathbf{z} \mid \boldsymbol{\lambda}, \varepsilon) = \sup_{\boldsymbol{\lambda}' \in B(\boldsymbol{\lambda}\mid\varepsilon)} \eta_k(\mathbf{z} \mid \boldsymbol{\lambda}')$. Lemma A.2 implies

$$\lim_{\varepsilon \to 0} \sum_k \overline{r}_k \cdot E\{\eta_k(\mathbf{Z}_k \mid \boldsymbol{\lambda}, \varepsilon)\} = \sum_k \overline{r}_k \cdot E\{\eta_k(\mathbf{Z}_k \mid \boldsymbol{\lambda})\}$$

and for any $\boldsymbol{\lambda} \in A$ Lemma A.1 gives

$$\sum_k \overline{r}_k \cdot E\{\eta_k(\mathbf{Z}_k \mid \boldsymbol{\lambda})\} < \sum_k \overline{r}_k \cdot E\{\eta_k(\mathbf{Z}_k \mid \boldsymbol{\lambda}^\circ)\}.$$

Hence there exists an $\varepsilon_{\boldsymbol{\lambda}} > 0$ such that

(A.15) $$\sum_k \overline{r}_k \cdot E\{\eta_k(\mathbf{Z}_k \mid \boldsymbol{\lambda}, \varepsilon_{\boldsymbol{\lambda}})\} < \sum_k \overline{r}_k \cdot E\{\eta_k(\mathbf{Z}_k \mid \boldsymbol{\lambda}^\circ)\}.$$

Since $A$ is compact, there are finitely many $\boldsymbol{\lambda}_1, \ldots, \boldsymbol{\lambda}_M \in A$ such that for any $\boldsymbol{\lambda} \in A$ there exists $1 \leq m \leq M$ with $\boldsymbol{\lambda} \in B^\circ(\boldsymbol{\lambda}_m \mid \varepsilon_{\boldsymbol{\lambda}_m})$. Thus $\eta_k(\mathbf{z} \mid \boldsymbol{\lambda}) \leq \eta_k(\mathbf{z} \mid \boldsymbol{\lambda}_m, \varepsilon_{\boldsymbol{\lambda}_m})$ and

(A.16) $$\sup_{\boldsymbol{\lambda} \in A} \ell^{(n)}(\boldsymbol{\lambda}) - \ell^{(n)}(\boldsymbol{\lambda}^\circ) \leq \operatorname*{Max}_m \sum_k \sum_i [\eta_k(\mathbf{Z}_{ki} \mid \boldsymbol{\lambda}_m, \varepsilon_{\boldsymbol{\lambda}_m}) - \eta_k(\mathbf{Z}_{ki} \mid \boldsymbol{\lambda}^\circ)].$$

For each $m$ the strong law of large numbers gives almost surely

$$\lim_{n \to 0} \frac{1}{n} \sum_k \sum_i [\eta_k(\mathbf{Z}_{ki} \mid \boldsymbol{\lambda}_m, \varepsilon_{\boldsymbol{\lambda}_m}) - \eta_k(\mathbf{Z}_{ki} \mid \boldsymbol{\lambda}^\circ)]$$

$$= \sum_k \overline{r}_k \cdot [E\{\eta_k(\mathbf{Z}_k \mid \boldsymbol{\lambda}_m, \varepsilon_{\boldsymbol{\lambda}_m})\} - E\{\eta_k(\mathbf{Z}_k \mid \boldsymbol{\lambda}^\circ)\}] < 0 \qquad [\text{cf. (A.15)}]$$

with finite expectations by (A.10) and (A.14). Hence

$$\lim_{n \to 0} \sum_k \sum_i [\eta_k(\mathbf{Z}_{ki} \mid \boldsymbol{\lambda}_m, \varepsilon \boldsymbol{\lambda}_m) - \eta_k(\mathbf{Z}_{ki} \mid \boldsymbol{\lambda}^\circ)] = -\infty$$

and the right-hand side in (A.16) tends to $-\infty$ for $n \to \infty$ almost surely. □

PROOF OF THEOREM 4 (CONSISTENCY). For any $\varepsilon > 0$, the function $\ell^{(n)}(\boldsymbol{\lambda})$ attains its maximum within $B(\boldsymbol{\lambda}^\circ \mid \varepsilon)$. We show first that (almost surely) the maximizing argument lies (for almost all $n$) in the *open* ball $B^\circ(\boldsymbol{\lambda} \mid \varepsilon)$, and hence is a solution of $D_{\boldsymbol{\lambda}} \ell^{(n)}(\boldsymbol{\lambda}) = \mathbf{0}$. Applying Lemma A.3 to the boundary $A_\varepsilon = \partial B(\boldsymbol{\lambda}^\circ \mid \varepsilon)$ yields that the following statements hold almost surely for almost all $n$:

(i) $\sup_{\boldsymbol{\lambda} \in A_\varepsilon} \ell^{(n)}(\boldsymbol{\lambda}) < \ell^{(n)}(\boldsymbol{\lambda}^\circ)$,
(ii) $\sup_{\|\boldsymbol{\lambda}-\boldsymbol{\lambda}^\circ\| \leq \varepsilon} \ell^{(n)}(\boldsymbol{\lambda}) \leq \sup_{\|\boldsymbol{\lambda}-\boldsymbol{\lambda}^\circ\| < \varepsilon} \ell^{(n)}(\boldsymbol{\lambda})$,
(iii) there exists $\tilde{\boldsymbol{\lambda}}^{(n)} \in B^\circ(\boldsymbol{\lambda}^\circ \mid \varepsilon)$ with $D_{\boldsymbol{\lambda}} \ell^{(n)}(\tilde{\boldsymbol{\lambda}}^{(n)}) = \mathbf{0}$,
(iv) $\ell^{(n)}(\boldsymbol{\lambda})$ is strictly concave [cf. (6.1), (R2′)],
(v) there is a unique $\tilde{\boldsymbol{\lambda}}^{(n)} \in B^\circ(\boldsymbol{\lambda}^\circ \mid \varepsilon)$ maximizing $\ell^{(n)}(\boldsymbol{\lambda})$,
(vi) $\tilde{\boldsymbol{\lambda}}^{(n)} = \hat{\boldsymbol{\lambda}}^{(n)}$ [cf. (C1)].

This proves (a), (b) and also (c), since $\varepsilon$ was arbitrary. □

### A.4. Proof of the results in Section 7.

PROOF OF THEOREM 5 (NORMALITY). The (standard) proof is only outlined. For $\mathbf{U}^{(n)} = D_{\boldsymbol{\lambda}}^T \ell^{(n)}$ the central limit theorem and (5.6) give

(A.17) $$n^{-1/2} \mathbf{U}^{(n)}(\boldsymbol{\lambda}^\circ) \xrightarrow[n \to \infty]{\mathcal{L}} N(\mathbf{0}, \overline{\boldsymbol{\Sigma}}(\boldsymbol{\lambda}^\circ)).$$

A first-order expansion about $\boldsymbol{\lambda}^\circ$ yields

$$n^{-1/2} \mathbf{U}^{(n)}(\hat{\boldsymbol{\lambda}}^{(n)}) = n^{-1/2} \mathbf{U}^{(n)}(\boldsymbol{\lambda}^\circ) + \mathbf{D}_n \cdot \sqrt{n}[\hat{\boldsymbol{\lambda}}^{(n)} - \boldsymbol{\lambda}^\circ]$$

with

$$\mathbf{D}_n := \frac{1}{n} \int_0^1 D_{\boldsymbol{\lambda}} \mathbf{U}^{(n)}(\boldsymbol{\lambda}^\circ + t[\hat{\boldsymbol{\lambda}}^{(n)} - \boldsymbol{\lambda}^\circ])\, dt$$

$$= -\frac{1}{n} \int_0^1 \mathbf{J}^{(n)}(\boldsymbol{\lambda}^\circ + t[\hat{\boldsymbol{\lambda}}^{(n)} - \boldsymbol{\lambda}^\circ])\, dt$$

and Condition (N1) implies

$$\mathbf{D}_n \cdot \sqrt{n}[\hat{\boldsymbol{\lambda}}^{(n)} - \boldsymbol{\lambda}^\circ] + n^{-1/2} \mathbf{U}^{(n)}(\boldsymbol{\lambda}^\circ) \xrightarrow[n \to \infty]{P} \mathbf{0}.$$

$\mathbf{D}_n$ can be replaced by its limit $-\overline{\mathbf{I}}(\boldsymbol{\lambda}^\circ)$ from Condition (N3), that is,

$$\sqrt{n}[\hat{\boldsymbol{\lambda}}^{(n)} - \boldsymbol{\lambda}^\circ] - n^{-1/2} \overline{\mathbf{I}}^{-1}(\boldsymbol{\lambda}^\circ) \mathbf{U}^{(n)}(\boldsymbol{\lambda}^\circ) \xrightarrow[n \to \infty]{P} \mathbf{0},$$

which together with (A.17) establishes (a). And (b) follows in view of Theorem 3(b). □

PROOF OF THEOREM 6. Keeping the notation from Appendix A.3, the partial derivatives of

$$\eta_k(\mathbf{z} \mid \boldsymbol{\lambda}) = \log \tilde{p}_k(\mathbf{z} \mid \boldsymbol{\lambda}) = \gamma_k^* + G_k(\mathbf{z}, \boldsymbol{\theta}) - \log\left(\sum_l \exp[\gamma_l^* + G_l(\mathbf{z}, \boldsymbol{\theta})]\right)$$

up to order 3 are given by

$$D_{\gamma_m^*}\eta_k(\mathbf{z} \mid \boldsymbol{\lambda}) = \Delta_{km} - \tilde{p}_m(\mathbf{z} \mid \boldsymbol{\lambda}),$$

$$D^2_{\gamma_m^* \lambda_s}\eta_k(\mathbf{z} \mid \boldsymbol{\lambda}) = -D_{\lambda_s}\tilde{p}_k(\mathbf{z} \mid \boldsymbol{\lambda}),$$

$$D^3_{\gamma_m^* \lambda_s \lambda_t}\eta_k(\mathbf{z} \mid \boldsymbol{\lambda}) = -D^2_{\lambda_s \lambda_t}\tilde{p}_k(\mathbf{z} \mid \boldsymbol{\lambda}),$$

$$D_{\theta_r}\eta_k(\mathbf{z} \mid \boldsymbol{\lambda}) = \sum_l [\Delta_{kl} - \tilde{p}_l(\mathbf{z} \mid \boldsymbol{\lambda})] \cdot D_{\theta_r}G_l(\mathbf{z}, \boldsymbol{\theta}),$$

$$D^2_{\theta_r \theta_s}\eta_k(\mathbf{z} \mid \boldsymbol{\lambda}) = \sum_l ([\Delta_{kl} - \tilde{p}_l(\mathbf{z} \mid \boldsymbol{\lambda})] \cdot D^2_{\theta_r \theta_s}G_l(\mathbf{z}, \boldsymbol{\theta})$$
$$- D_{\theta_s}\tilde{p}_l(\mathbf{z} \mid \boldsymbol{\lambda}) \cdot D_{\theta_r}G_l(\mathbf{z}, \boldsymbol{\theta})),$$

$$D^3_{\theta_r \theta_s \theta_t}\eta_k(\mathbf{z} \mid \boldsymbol{\lambda}) = \sum_l ([\Delta_{kl} - \tilde{p}_l(\mathbf{z} \mid \boldsymbol{\lambda})] \cdot D^3_{\theta_r \theta_s \theta_t}G_l(\mathbf{z}, \boldsymbol{\theta})$$
$$- D_{\theta_t}\tilde{p}_l(\mathbf{z} \mid \boldsymbol{\lambda}) \cdot D^2_{\theta_r \theta_s}G_l(\mathbf{z}, \boldsymbol{\theta})$$
$$- D_{\theta_s}\tilde{p}_l(\mathbf{z} \mid \boldsymbol{\lambda}) \cdot D^2_{\theta_r \theta_t}G_l(\mathbf{z}, \boldsymbol{\theta})$$
$$- D^2_{\theta_s \theta_t}\tilde{p}_l(\mathbf{z} \mid \boldsymbol{\lambda}) \cdot D_{\theta_r}G_l(\mathbf{z}, \boldsymbol{\theta}))$$

with partial derivatives [cf. (A.1)]

$$D_{\lambda_s}\tilde{p}_k(\mathbf{z} \mid \boldsymbol{\lambda}) = \tilde{p}_k(\mathbf{z} \mid \boldsymbol{\lambda}) \cdot D_{\lambda_s}\eta_k(\mathbf{z} \mid \boldsymbol{\lambda}),$$
$$D^2_{\lambda_s \lambda_t}\tilde{p}_k(\mathbf{z} \mid \boldsymbol{\lambda}) = \tilde{p}_k(\mathbf{z} \mid \boldsymbol{\lambda})[D^2_{\lambda_s \lambda_t}\eta_k(\mathbf{z} \mid \boldsymbol{\lambda}) + D_{\lambda_s}\eta_k(\mathbf{z} \mid \boldsymbol{\lambda}) \cdot D_{\lambda_t}\eta_k(\mathbf{z} \mid \boldsymbol{\lambda})].$$

Next we deduce from Condition (MC) a weaker moment condition, from which Conditions (N3) and (N4) will be derived (cf. Lemma A.4):

CONDITION (MC)~. There exists $\varepsilon^\circ > 0$ such that for $B(\boldsymbol{\lambda}^\circ) = \{\boldsymbol{\lambda} \mid \|\boldsymbol{\lambda} - \boldsymbol{\lambda}^\circ\| \leq \varepsilon^\circ\}$ and all $k = 0, \ldots, K$ the following functions:

$$\sup_{\boldsymbol{\lambda} \in B(\boldsymbol{\lambda}^\circ)} |D^3_{\lambda_r \lambda_s \lambda_t}\eta_l(\mathbf{Z}_k \mid \boldsymbol{\lambda})| \qquad \text{with } \mathbf{Z}_k = h_X(X_k)$$

have finite expectation for all $r, s, t = 1, \ldots, S$ and $l = 0, \ldots, K$.

For the above derivatives we successively get the following bounds, where the fixed argument $\mathbf{z}$ is omitted:

$$|D_{\gamma_m^*}\eta_k(\boldsymbol{\lambda})| \leq 1, \qquad |D_{\theta_r}\eta_k(\boldsymbol{\lambda})| \leq H_+(\boldsymbol{\theta}),$$

$$|D_{\lambda_r}\eta_k(\boldsymbol{\lambda})| \leq H_+^*(\boldsymbol{\theta}) := 1 + H_+(\boldsymbol{\theta}),$$

$$|D^2_{\gamma_m^*\lambda_s}\eta_k(\boldsymbol{\lambda})| = |D_{\lambda_s}\tilde{p}_k(\boldsymbol{\lambda})| \leq |D_{\lambda_s}\eta_k(\boldsymbol{\lambda})| \leq H_+^*(\boldsymbol{\theta}),$$

$$|D^2_{\theta_r\theta_s}\eta_k(\boldsymbol{\lambda})| \leq H_{++}(\boldsymbol{\theta}) + H_+(\boldsymbol{\theta})^2,$$

$$|D^2_{\lambda_s\lambda_t}\eta_k(\boldsymbol{\lambda})| \leq H_+^*(\boldsymbol{\theta})^2 + H_{++}(\boldsymbol{\theta}),$$

$$|D^2_{\lambda_s\lambda_t}\tilde{p}_k(\boldsymbol{\lambda})| \leq 2H_+^*(\boldsymbol{\theta})^2 + H_{++}(\boldsymbol{\theta}),$$

$$|D^3_{\gamma_m^*\lambda_s\lambda_t}\eta_k(\boldsymbol{\lambda})| = |D^2_{\lambda_s\lambda_t}\tilde{p}_k(\boldsymbol{\lambda})| \leq 2H_+^*(\boldsymbol{\theta})^2 + H_{++}(\boldsymbol{\theta}),$$

$$|D^3_{\theta_r\theta_s\theta_t}\eta_k(\boldsymbol{\lambda})| \leq H_{+++}(\boldsymbol{\theta}) + 3H_+^*(\boldsymbol{\theta})H_{++}(\boldsymbol{\theta}) + 2H_+^*(\boldsymbol{\theta})^3.$$

Taking (for fixed $\mathbf{z}$) the supremum over the ball $B(\boldsymbol{\theta}^\circ)$ gives

$$\sup H_+^* \leq 1 + \sum_r \sup H_r,$$

$$\sup H_+^{*2} \leq 1 + 2\sum_s \sup H_s + \sum_s\sum_t \sup H_{st},$$

$$\sup H_+^{*3} \leq 1 + 3\sum_r \sup H_r + 3\sum_r\sum_s \sup H_r H_s + \sum_r\sum_s\sum_t \sup H_r H_s H_t,$$

$$\sup H_{++} \leq \sum_s\sum_t \sup H_{st},$$

$$\sup H_{+++} \leq \sum_r\sum_s\sum_t \sup H_{rst},$$

$$\sup H_+^* \cdot H_{++} \leq \sum_r\sum_s\sum_t [\sup H_{st} + \sup H_r H_{st}].$$

Condition (MC) obviously implies for $i = 1, 2$ that

$$\sup_{\boldsymbol{\theta}\in B(\boldsymbol{\theta}^\circ)} H_r(\mathbf{Z}_k \mid \boldsymbol{\theta})^i, \qquad \sup_{\boldsymbol{\theta}\in B(\boldsymbol{\theta}^\circ)} H_r(\mathbf{Z}_k \mid \boldsymbol{\theta}) \cdot H_{st}(\mathbf{Z}_k \mid \boldsymbol{\theta})$$

have finite expectation, too. Hence

$$\sup_{\boldsymbol{\gamma}} \sup_{\boldsymbol{\theta}\in B(\boldsymbol{\theta}^\circ)} |D^3_{\lambda_r\lambda_s\lambda_t}\eta_l(\mathbf{Z}_k \mid \boldsymbol{\theta}, \boldsymbol{\gamma})|$$

has finite expectation for any $r, s, t$ and any $k, l$. This proves Condition (MC)$^\sim$ and Lemma A.4 establishes the theorem. $\square$

LEMMA A.4. *Conditions (N2) and (MC)$^\sim$ imply Conditions (N3) and (N4).*

PROOF. Using (6.1) for $\boldsymbol{\lambda} = \boldsymbol{\lambda}^\circ$ to establish Condition (N3), it suffices to show for any $s$ and $t$ that

$$(A.18) \quad \frac{1}{n} \int_0^1 [J_{st}^{(n)}(\boldsymbol{\lambda}^\circ + t[\hat{\boldsymbol{\lambda}}^{(n)} - \boldsymbol{\lambda}^\circ]) - J_{st}^{(n)}(\boldsymbol{\lambda}^\circ)]\, dt \xrightarrow[n\to\infty]{P} 0.$$

From

$$J_{st}^{(n)}(\boldsymbol{\lambda}) = -\sum_k \sum_i D^2_{\lambda_s \lambda_t} \eta_k(\mathbf{Z}_{ki} \mid \boldsymbol{\lambda}) \quad \text{with } \mathbf{Z}_{ki} = h_X(X_{ki})$$

a Taylor expansion gives for any $\varepsilon > 0$ and $\|\boldsymbol{\lambda} - \boldsymbol{\lambda}^\circ\| < \varepsilon$

$$\frac{1}{n} |J_{st}^{(n)}(\boldsymbol{\lambda}) - J_{st}(\boldsymbol{\lambda}^\circ)| \leq \varepsilon \overline{S}^{(n)}(\varepsilon),$$

$$\overline{S}^{(n)}(\varepsilon) = \frac{1}{n} \sum_k \sum_i \sup_{\|\boldsymbol{\lambda}' - \boldsymbol{\lambda}^\circ\| \leq \varepsilon} \|D^3_{\boldsymbol{\lambda}\lambda_s \lambda_t} \eta_k(\mathbf{Z}_{ki} \mid \boldsymbol{\lambda}')\|.$$

The strong law of large numbers yields

$$\overline{S}^{(n)}(\varepsilon) \xrightarrow[n\to\infty]{} \sum_k \overline{r}_k E\left(\sup_{\|\boldsymbol{\lambda}' - \boldsymbol{\lambda}^\circ\| \leq \varepsilon} \|D^3_{\boldsymbol{\lambda}\lambda_s \lambda_t} \eta_k(\mathbf{Z}_k \mid \boldsymbol{\lambda}')\|\right) \quad \text{almost surely,}$$

where the limit is finite by Condition (MC)$^\sim$ for $\varepsilon \leq \varepsilon^\circ$. For $\|\hat{\boldsymbol{\lambda}}^{(n)} - \boldsymbol{\lambda}^\circ\| < \varepsilon$ we thus have

$$\left| \frac{1}{n} \int_0^1 [J_{st}^{(n)}(\boldsymbol{\lambda}^\circ + t[\hat{\boldsymbol{\lambda}}^{(n)} - \boldsymbol{\lambda}^\circ]) - J_{st}^{(n)}(\boldsymbol{\lambda}^\circ)]\, dt \right| \leq \frac{1}{n} \sup_{\|\boldsymbol{\lambda} - \boldsymbol{\lambda}^\circ\| \leq \varepsilon} |J_{st}^{(n)}(\boldsymbol{\lambda}) - J_{st}(\boldsymbol{\lambda}^\circ)|$$

$$\leq \varepsilon \overline{S}^{(n)}(\varepsilon)$$

which in view of Condition (N2) implies (A.18). And Condition (N4) follows similarly. Note that if almost sure convergence $\hat{\boldsymbol{\lambda}}^{(n)} \to \boldsymbol{\lambda}^\circ$ is assumed instead of Condition (N2), then the above arguments establish almost sure convergence in Conditions (N3) and (N4), too. $\square$

**Acknowledgments.** The author thanks both referees for their helpful comments which improved the first draft of the paper.

## REFERENCES

[1] BRESLOW, N. E., ROBINS, J. M. and WELLNER, J. A. (2000). On the semiparametric efficiency of logistic regression under case-control sampling. *Bernoulli* **6** 447–455. MR1762555
[2] CSISZÁR, I. (1975). *I*-divergence geometry of probability distributions and minimization problems. *Ann. Probab.* **3** 146–158. MR0365798
[3] GILULA, Z. and HABERMAN, S. J. (1988). The analysis of multivariate contingency tables by restricted canonical and restricted association models. *J. Amer. Statist. Assoc.* **83** 760–771. MR0963804

[4] GOODMAN, L. A. (1985). The analysis of cross-classified data having ordered and/or unordered categories: Association models, correlation models, and asymmetry models for contingency tables with or without missing entries. *Ann. Statist.* **13** 10–69. MR0773152
[5] HABERMAN, S. J. (1974). *The Analysis of Frequency Data*. Univ. Chicago Press, Chicago. MR0408098
[6] KAGAN, A. (1988). A note on the logistic function. *Biometrika* **88** 599–601. MR1844857
[7] KULLBACK, S. (1959). *Information Theory and Statistics.* Dover, Mineola, NY. MR1461541
[8] MCCULLAGH, P. and NELDER, J. A. (1989). *Generalized Linear Models*, 2nd ed. Chapman and Hall, London. MR0727836
[9] OSIUS, G. (2000). The association between two random elements: A complete characterization in terms of odds ratios. Mathematik-Arbeitspapiere 53, Univ. Bremen. Available at http://www.math.uni-bremen.de/~osius.
[10] OSIUS, G. (2004). The association between two random elements: A complete characterization and odds ratio models. *Metrika* **60** 261–277. MR2189755
[11] PLACKETT, R. L. (1974). *The Analysis of Categorical Data.* Griffin, London.
[12] PRENTICE, R. L. and PYKE, R. (1979). Logistic disease incidence models and case-control studies. *Biometrika* **66** 403–411. MR0556730
[13] RÜSCHENDORF, L. and THOMSEN, W. (1993). Note on the Schrödinger equation and I-projections. *Statist. Probab. Lett.* **17** 369–375. MR1237783
[14] SCOTT, A. J. and WILD, C. J. (1997). Fitting logistic regression models in stratified case-control studies. *Biometrika* **84** 57–71. MR1450191
[15] WALD, A. (1949). Note on the consistency of the maximum likelihood estimate. *Ann. Statist.* **20** 595–601. MR0032169
[16] WEINBERG, C. R. and WACHOLDER, S. (1993). Prospective analysis of case-control data under general multiplicative-intercept risk models. *Biometrika* **80** 461–465. MR1243520

INSTITUT FÜR STATISTIK
FACHBEREICH 3
UNIVERSITÄT BREMEN
GERMANY
E-MAIL: osius@math.uni-bremen.de